%% file: Paperling_revu.tex
\documentclass[unknownkeysallowed,11pt]{article}
\usepackage[final]{graphicx}
\usepackage[utf8]{inputenc}
\usepackage[english]{babel}
\usepackage[T1]{fontenc}
\input{macropulko}

\usepackage{framed}

\usepackage{calc}
\usepackage[usenames, dvipsnames]{xcolor}
\usepackage{wasysym}
\usepackage{pifont}
\usepackage{upgreek}
\usepackage{xparse}
\usepackage{fourier-orns} 
\usepackage{tikz, ifthen}
	\usetikzlibrary{arrows,shapes,positioning,automata,decorations.markings}
\usepackage{float}
\tikzstyle directed=[postaction={decorate,decoration={markings,
    mark=at position .65 with {\arrow{latex}}}}]
\usepackage[margin=20pt,format=plain,small,justification=centering]{caption}
\usepackage{moreverb}
\usepackage{amsmath,amscd}
\usepackage{amsthm}
\usepackage{amsfonts}
\usepackage{amssymb}
\usepackage{pgfplots}

\usepackage{fourier} 

\usepackage{array}
\usepackage{bmpsize}
\usepackage{mdwlist} 
\usepackage{wrapfig}
\usepackage{shadethm}
\usepackage{setspace}
\usetikzlibrary{arrows,%
                petri,%
                topaths}%
\usepackage{tkz-graph}
\usepackage{tkz-berge}

\usepackage[left=2.4 cm,right=2.4cm,top=2.2cm,bottom=2.2cm]{geometry}
\usepackage{multicol}
\usepackage{fancyhdr}
\usepackage{mathtools}

\newtheoremstyle{thm}
  {12pt}
  {12pt}
  {\slshape}
  {}
  {\bfseries \scshape}
  {. ---}
  { }
  {}

\NewDocumentEnvironment{hypo}{o}
  {\definecolor{shadethmcolor}{rgb}{255,255,255} 
\definecolor{shaderulecolor}{rgb}{0,0.0,0.0} 
\setlength{\shadeboxrule}{1pt} 
    \IfNoValueTF{#1}
      {\begin{hyp}}
      {\begin{hyp}\textup{[#1]}}
  }
  {\end{hyp}\vspace{2pt}
  }
  
  \NewDocumentEnvironment{addass}{o}
  {\definecolor{shadethmcolor}{rgb}{255,255,255} 
\definecolor{shaderulecolor}{rgb}{0,0.0,0.0} 
\setlength{\shadeboxrule}{1pt} 
    \IfNoValueTF{#1}
      {\begin{addas}}
      {\begin{addas}\textup{[#1]}}
  }
  {\end{addas}\vspace{2pt}
  }

  \NewDocumentEnvironment{defmodd}{o}
  {\definecolor{shadethmcolor}{rgb}{255,255,255} 
\definecolor{shaderulecolor}{rgb}{0,0.0,0.0} 
\setlength{\shadeboxrule}{1pt} 
    \IfNoValueTF{#1}
      {\begin{defmod}}
      {\begin{defmod}\textup{[#1]}}
  }
  {\end{defmod}\vspace{2pt}
  }

\theoremstyle{thm} 
\newshadetheorem{theorem}{Theorem}
\newshadetheorem{hyp}{Assumption}

\newtheorem{addas}{Assumption}

\newtheorem*{defmod}{Definition}

\newshadetheorem{methodegrise}[theorem]{Méthode}
\newshadetheorem{defpropgrise}[theorem]{Definition-proposition}
\newshadetheorem{propgrise}[theorem]{Proposition}
\newtheorem{prop}[theorem]{Proposition}

\newtheorem{defin}[theorem]{Definition}

\newtheoremstyle{rq}{}{}{}{}{\bfseries \scshape}{.}{.5em}{}
\theoremstyle{rq}

\newcounter{ex}
\newtheorem{exemple}[ex]{Example}

\newtheoremstyle{exerci}{}{}{}{}{\bfseries}{.}{.5em}{}
\theoremstyle{exerci}

\newtheoremstyle{solu}{}{}{}{}{\bfseries}{.}{.5em}{#1\thmnote{#3}}
\theoremstyle{solu}
\newtheorem{demprop}{Proof of Proposition  }

\newtheorem{demthm}{Proof of Theorem  }

\newtheoremstyle{prob}{2cm}{}{}{6.32cm}{\scshape \bfseries}{.\newline\newline}{.5em}{}
\theoremstyle{prob}

\AtBeginDocument{}
\usepackage{enumitem}
\usepackage{xparse}
\let\CLASSsection\section
\RenewDocumentCommand{\section}{som}{%
  \IfBooleanTF{#1}
   {
    \CLASSsection*{#3}\markboth{{\textup{\textsc{#3}}}}{}%
   }
   {
    \IfNoValueTF{#2}
     {
      \CLASSsection{#3}%
     }
     {
      \CLASSsection[#2]{#3}%
     }%
   }%
}
\makeatletter
\renewcommand{\tableofcontents}{%
  \section*{\contentsname}%
  \@starttoc{toc}%
}
\makeatother

\renewcommand{\sectionmark}[1]{%
\markboth{\thesection.\ \textup{\textsc{#1}}}{}}

\usepackage{fnpct}
\usepackage{titling}
\thanksmarkseries{arabic}

\title{Stability of Piecewise Deterministic Markovian \\ Metapopulation Processes on Networks}
\author{Pierre Montagnon\thanks{~CMAP, École Polytechnique, Route de Saclay, 91128 Palaiseau Cedex, France.}~\textsuperscript{ ,}\thanks{~MaIAGE, INRA, Université Paris-Saclay, 78350 Jouy-en-Josas, France.}}
\begin{document}

\maketitle

\begin{abstract}
The purpose of this paper is to study a Markovian metapopulation model on a directed graph with edge-supported transfers and deterministic intra-nodal population dynamics. We first state tractable stability conditions for two typical frameworks motivated by applications: constant jump rates with multiplicative transfer amplitudes, and coercive jump rates with unitary transfers. More general criteria for boundedness, petiteness and ergodicity are then given.
\end{abstract}

{\bf Keywords:} metapopulation models, population dynamics, piecewise deterministic Markov processes, stability conditions, ergodicity

\sectionmark{Introduction}
\section{Introduction}

\indent \indent Metapopulation models describe the behavior of a population naturally or artificially split into spatially distant patches connected through individual movements \cite{Lev,Rit,MAW,Verb}. Although most such models have been dealing with fully deterministic or fully stochastic population dynamics, network-organized systems with deterministic intra-nodal dynamics and stochastic inter-nodal transfers naturally appear as relevant to describe metapopulations with low local stochasticity (e.g. because of large population sizes) and random population transfers between patches. 

\begin{center}
\begin{figure}[!h]
\begin{minipage}[c]{.35\linewidth}
\begin{tikzpicture}[xscale=.60,yscale=.75]
\tikzstyle{VertexStyle}=[scale=.75, shape        = circle,
                             fill         = red,
                             minimum size = 28pt,
                             text         = black,
                             draw]
  \Vertex[L=$2$, a=50 , d=3.5 cm]{A}
  \tikzstyle{VertexStyle}=[scale=.75, shape        = circle,
                             fill         = lime,
                             minimum size = 28pt,
                             text         = black,
                             draw]
  \Vertex[L=$3$, a=54 , d=6.5cm]{B} 
  \tikzstyle{VertexStyle}=[scale=.9, shape        = circle,
                             fill         = yellow!30!red!90,
                             minimum size = 35pt,
                             text         = black,
                             draw]
  \Vertex[L=$4$, a=83 , d=4.5 cm]{C} 
    \Vertex[L=$7$, a=140 , d=4.7 cm]{F} 
      \tikzstyle{VertexStyle}=[scale=.9, shape        = circle,
                             fill         = yellow!10!red!40,
                             minimum size = 35pt,
                             text         = black,
                             draw]
  \Vertex[L=$8$, a=230 , d=3.3 cm]{G} 
  \Vertex[L=$5$, a=108 , d=4.6 cm]{D} 

 \tikzstyle{VertexStyle}=[scale=.60, shape        = circle,
                             fill         = red,
                             minimum size = 28pt,
                             text         = black,
                             draw]
  \Vertex[L=$6$, a=124 , d=3.1 cm]{E} 
   \tikzstyle{VertexStyle}=[scale=.60, shape        = circle,
                             fill         = yellow!10!red!95,
                             minimum size = 28pt,
                             text         = black,
                             draw]
  \Vertex[L=$9$, a=252 , d=4.8 cm]{H}  
   \tikzstyle{VertexStyle}=[scale=.60, shape        = circle,
                             fill         = red!60!yellow!70,
                             minimum size = 28pt,
                             text         = black,
                             draw]
  \Vertex[L=$10$, a=288 , d=4.4 cm]{I}  
  
  \tikzstyle{VertexStyle}=[scale=.90,shape        = circle,
                             fill         = red,
                             minimum size = 20pt,
                             text         = black,
                             draw]
                             \Vertex[L=$11$, a=324 , d=4.5 cm]{J}  
 
      \tikzstyle{VertexStyle}=[shape        = circle,
                             fill         = lime,
                             minimum size = 50pt,
                             text         = black,
                             draw]
  \Vertex[L=$1$, a=0 , d=0 cm]{Z}  
            \tikzstyle{EdgeStyle}=[post,>=latex,->,very thick, bend right]    
        \Edge(Z)(D)
        \tikzstyle{EdgeStyle}=[post,>=latex,->,very thick, bend left]
          \Edge(Z)(F)
                \Edge(Z)(H)
                  \tikzstyle{EdgeStyle}=[post,>=latex,->,ultra thick, bend left]
                           \Edge(Z)(I)
                                 \tikzstyle{EdgeStyle}=[post,>=latex,->,ultra thick, bend right]
                  \Edge(Z)(B)
                    \tikzstyle{EdgeStyle}=[post,>=latex,->,thick, bend left]
                    \Edge(Z)(J)
                        \tikzstyle{EdgeStyle}=[post,>=latex,->,thin, bend right]
                          \Edge(A)(C)
    \Edge(B)(C)
      \Edge(C)(A)
                      \Edge(F)(G)
                       \tikzstyle{EdgeStyle}=[post,>=latex,->,ultra thin, bend right]
    \Edge(D)(E)
              \Edge(G)(H)
\end{tikzpicture}\\\begin{center}
{\huge A}
\end{center}\vspace{1mm}   \end{minipage}
\begin{minipage}[c]{.55\linewidth}
\includegraphics[scale=.20]{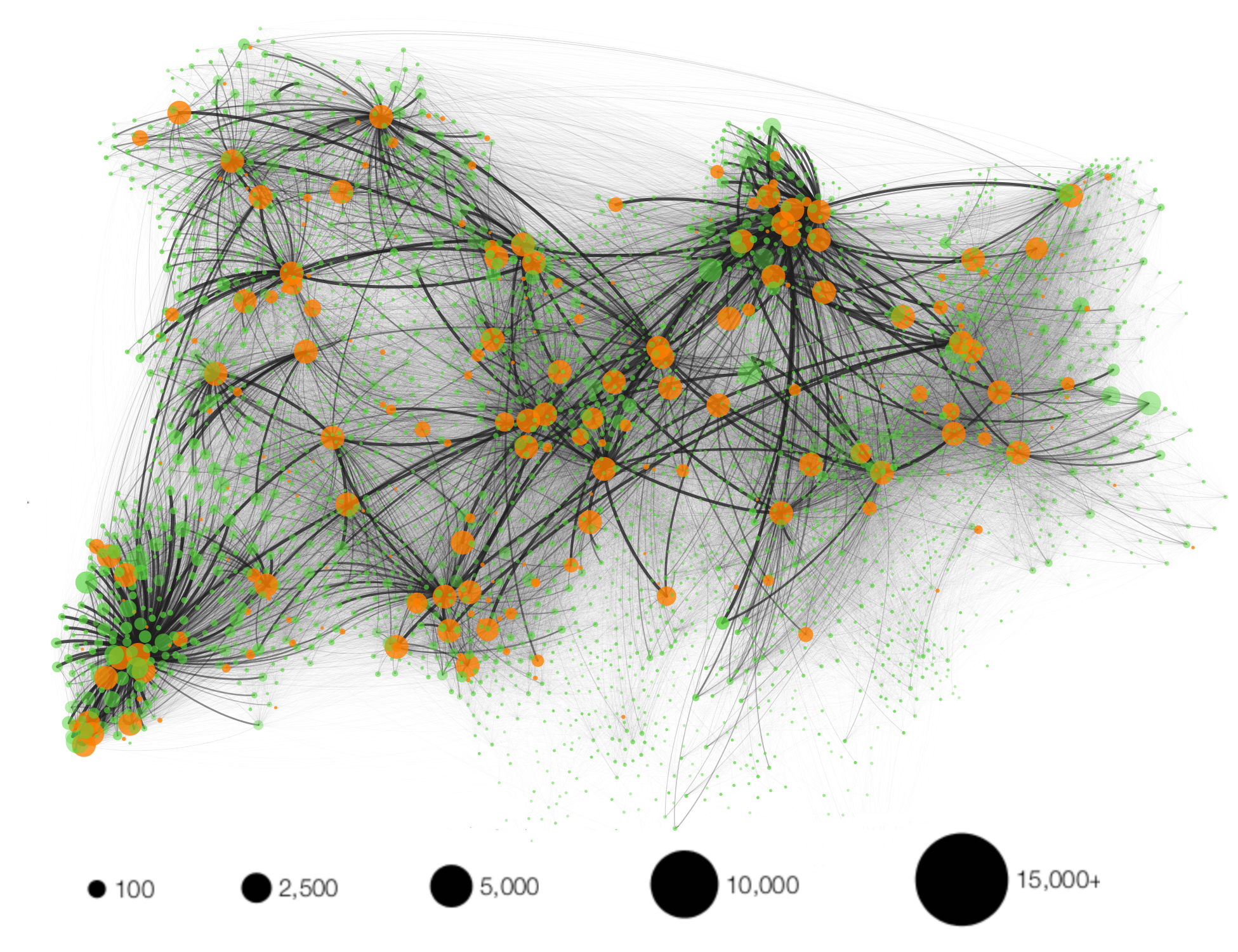}\\\begin{center}
{\huge B}
\end{center}\vspace{1mm}
   \end{minipage}
\captionsetup{justification=justified}
\caption{\label{fig_illus_intro}Panel A: stochastic "mainland-islands" Levins model \cite{Lev} with unitary transfers. Edge thickness corresponds to transfer rates, while colors denote the sign of autonomous growth functions (red tones being associated to positive growth functions, green tones to positive ones). Panel B: 2008 cattle trade network in the French Auvergne-Rhône-Alpes area, densely population with cattle (courtesy of G. Beaunée). Edge thickness corresponds to the overall volume of cattle transfers within the period and colors stand for holding types (green for farms and orange for commercial operators).}
\end{figure} 
\end{center}

The original motivation of this paper is the study of a population model on a cattle trade network where nodes are farms and commercial operators (see \cite{Dut,Hos} and Figure \ref{fig_illus_intro}). Entries in such a system are deterministic and individuals move between nodes as they would in a Jackson network \cite{Jack,Mjack,Dai}. Similar models may allow for describing human transportation and animal movements causing epidemic propagation \cite{Col,Ker,NGT,BPR}.

\indent There is, to our knowledge, little literature on general stability criteria for $\R^n$-valued piecewise deterministic Markov processes. Meaningful results on the stability of Jackson networks have been derived that are deeply rooted on considerations about the graph structure \cite{Wal}. We wish to obtain similar statements in our semi-deterministic framework while allowing for state-dependent jump intensities and possibly large transfers. The total population of the system is here preserved by jumps and behaves as a non-Markovian randomly switched process \cite{Lib, Mal}. That is why we expect the process boundedness to arise from conditions on the deterministic inter-jump flow and the process ability to reach population-decreasing states in short times.

\indent \indent Our first task is to define a general metapopulation model on a network. We will represent the $\R_+^n$-valued population of the system patches (that is, nodes) as a piecewise deterministic Markov process $X=(X_1,\ldots,X_n)$. The population of each patch $i$ will be associated with an autonomous growth function $\phi_i$ (meaning that $\mathrm{d}X_i(t)=\phi_i(X(t))\mathrm{d}t$ between jumps), and an instantaneous transfer from patch $i$ to patch $j$ at population $x$ occurs at state-dependent rate $\theta_{i,j}(x)$, its amplitude being drawn according to a $[0,x_i]$-supported law $\mu_{i,j}(x,\cdot)$. We aim at providing sufficient conditions for the process stability as well as some asymptotic properties in the stable case.\pe

\indent \indent Modelling cattle trade network dynamics and trying to allow for macroscopic jumps led to the two following motivating examples. If $\mathcal{G}=(\ie{1}{n},\mathcal{A})$ is a strongly connected directed graph on $\ie{1}{n}$, define $(c_1,\ldots,c_n)\in \R_+^n\setminus\{0\}$ and consider the constant-growth settings defined by:
\begin{center}
\textbf{Multiplicative uniform setting}\vspace{-1mm}
\end{center}\begin{equation}\label{modele_motiv_multipl}\phi_i(x)=\begin{cases}c_i \text{ if }c_i\geqslant 0 \\ c_i\mathrm{1}_{x_i>0} \text{ if }c_i<0\end{cases} \quad \text{and} \quad \begin{cases} \theta_{i,j}(x)=\theta_{i,j} \\ \mu_{i,j}(x,\cdot)=\mathcal{U}([0,x_i])\end{cases}\end{equation} with $\theta_{i,j}>0$ if $(i,j)\in \mathcal{A}$ and $\theta_{i,j}= 0$ otherwise, or: 
\begin{center}
\textbf{Unitary power setting}\vspace{-1mm}
\end{center}
\begin{equation}\label{modele_motiv_unit}\phi_i(x)=\begin{cases}c_i \text{ if }c_i\geqslant 0 \\ c_i\mathrm{1}_{x_i>0} \text{ if }c_i<0\end{cases} \quad \text{and} \quad \begin{cases} \theta_{i,j}(x)=(1\vee x_i)^{\alpha} \\ \mu_{i,j}(x,\cdot)=\delta_{1\wedge x_i}\end{cases}\end{equation} for some $\alpha\in (0,1]$.

\indent \indent Specification (\ref{modele_motiv_multipl}) corresponds to the case of patches transferring at constant temporal rate a random fraction of their population to each other, which fits to empirical data on cattle trade observed over large time steps. Specification (\ref{modele_motiv_unit}) corresponds to unitary population transfers occurring at population-dependent rates, which is the natural framework for modelling cattle trade movements observed on a daily basis since the actual size of cattle transfers is bounded by transportation capacities. Note that the latter setting only differs from an open Jackson network by its continuous state space and the existence of the constant growth flow.

\indent \indent Theorems \ref{thm_pcpal_mult} and \ref{thm_pcpal_uni} below imply that the population process is positive Harris recurrent if and only if \begin{equation}\label{cond_flux_tot}\sum_{i=1}^n c_i<0\end{equation}
in which case it is also $F$-ergodic for some exponential function $F:\R_+^n\to \R_+$. Note that equation (\ref{cond_flux_tot}) is reminiscent of the usual traffic condition for queues networks stated in equation (1.9) of \cite{Dai} but does not involve any term related to inter-patch transfers. This stems from the fact that transfers along $\mathcal{G}$ happen and spread the total population among all patches, thus letting the flow on $\left(\R_+^*\right)^n$ bring $X$ back towards lower population states. We will see that such a mechanism prevails in more general settings whenever jumps are large or frequent enough; in this case, a simple condition on the autonomous growth functions is required to ensure the process stability.
\pe \\
\indent This paper is organized as follows. We first give a formal definition of our general piecewise deterministic metapopulation process (Section \ref{section_modelisation}) of which we highlight two instances of interest referred to as \emph{multiplicative} and \emph{unitary} frameworks (Section \ref{section_models_multuni}), and state stability results for these settings. We then give Meyn-Tweedie inspired criteria for boundedness (Section \ref{section_bornitude}), petiteness (Section \ref{section_petitesse}) and ergodicity (Section \ref{section_ergo}) in the general case. Results from Section \ref{section_general_crit} apply to prove the results of Section \ref{section_models_multuni}. All proofs are postponed to Section \ref{section_demos}.
\pe \\
\indent From now on, let $\mathcal{I}=\{(i,j)\in \ie{1}{n}^2\mid i\neq j\}$. For any $i\in \ie{1}{n}$, we denote by $e_i$ the $i$-th vector of the canonical basis of $\R^n$. For any $d\geqslant 1$ and any $x\in \R_+^d$ (respectively any $\R^d$-valued process $X=(X(t))_{t\geqslant 0}$) we write $x_i$ (resp. $X_i=(X_i(t))_{t\geqslant 0}$) for the $i$-th coordinate of $x$ (resp. of $X$). For any topological space $T=(E,\mathcal{T})$, we denote by $\mathcal{B}(T)=\mathcal{B}(\mathcal{T})$ the Borel $\sigma$-algebra on $T$. If $d\geqslant 1$, $\lambda_d$ will stand for the Lebesgue measure on the $d$-dimensional affine subsets of the $\R^K$ spaces (with $d\leqslant K$). Finally, we use the standard convention $\inf(\varnothing)=+\infty$.

\section{A general piecewise deterministic metapopulation model}
\label{section_modelisation}

\indent \indent For all $i\in \ie{1}{n}$, let $\phi_i:\R_+^n\to \R$ be a measurable and bounded function. Assume that $\phi=(\phi_1,\ldots,\phi_n)$ is such that the flow $\Phi$ associated to the vector field $\phi$ is well-defined on $\R_+^n$, that is, that for any $x\in \R_+^n$ there exists a unique continuous function $\Phi(x,\cdot):\R_+\to \R_+^n$ such that $\Phi(x,t)=x+\int_0^t \phi(x,u)\mathrm{d}u$ for all $t\geqslant 0$. This condition is for instance fulfilled if $\phi$ is $\mathcal{C}^1$.\pe 

We remind that $\mathcal{I}=\{(i,j)\in \ie{1}{n}^2 : i\neq j\}$. For all $(i,j)\in \mathcal{I}$, let $\theta_{i,j}:\R_+^n\to \R$ be a positive measurable, locally bounded function. We assume that there exists a directed graph $\mathcal{G}=(\ie{1}{n},\mathcal{E})$, with $\mathcal{E}\subset \mathcal{I}$, such that for all $(i,j)\in \mathcal{E}$, $\theta_{i,j}(x)>0$ whenever $x_i>0$. We call $\mathcal{G}$ the \emph{active graph} of the model.\pe 

For all $(i,j)\in \mathcal{I}$ and all $x\in \R_+^n$, finally define a probability measure $\mu_{i,j}(x,\cdot)$ on $[0,x_i]$, and for all $\xi\in [0,1]$, set $$q_{i,j}(x,\xi)=\inf~\{u\in [0,x_i] \mid \mu_{i,j}(x,[0,u])\geqslant \xi\}$$
and assume that the $q_{i,j}:\R_+^n\times [0,1]\to \R_+$ are Borel-measurable.\pe 

Our object of interest is the family of processes solution of the following stochastic differential equation:\begin{equation}\label{EDS}\fbox{$\displaystyle{\mathrm{d}X(t) = \phi(X(t))\mathrm{d}t + \sum_{(i,j)\in \mathcal{I}}\int_{\R_+\times [0,1]} q_{i,j}(X(t^-)},\xi)\left(e_j-e_i\right) \mathrm{1}_{z<\theta_{i,j}(X(t^-))}N_{i,j}(\mathrm{d}t,\mathrm{d}z,\mathrm{d}\xi)$}\end{equation}
where $\left(N_{i,j}\right)_{(i,j)\in \mathcal{I}}$ is a collection of independent homogeneous Poisson point measures on $\R_+\times \R_+\times [0,1]$ with intensity $\mathrm{d}t\mathrm{d}z\mathrm{d}\xi$ defined on $(\Omega,\mathcal{A},\P)$.\pe

If $x\in \R_+^n$, the existence of a strong solution $X^x=(X^x(t))_{t\geqslant 0}$ to (\ref{EDS}) with initial value $x$ follows from the explicit construction by Davis (\cite{Dav93} p.55). Pathwise uniqueness and uniqueness in law also hold since such solutions can naturally be expressed as deterministic functions of the atoms of the $N_{i,j}$ random measures. Moreover, we see at once that the $X^x$ processes are non-explosive and well-defined on $\R_+$ because the $\theta_{i,j}$ and $\phi_i$ are bounded on compacts and because jumps conserve the total population. By Theorem 25.5 of \cite{Dav93} (slightly adapted to allow for zero-amplitude jumps), the strong Markov property holds for $(X^x)_{x\in \R_+^n}$.\pe

We assume that there exists a Markov family $\left(X,\left(\P_x\right)_{x\in \R_+^n}\right)$ on $(\Omega,\mathcal{A})$ such that for any $x\in \R_+^n$, the law of $X$ under $\P_x$ is the law of $X^x$ under $\P$. We will mostly use this homogeneous and lighter notation.\pe 

We will denote by $(T_k)_{k\geqslant 1}$ and $(U_k)_{k\geqslant 1}$ the sequences of jump times and jump quantiles of $X$, the second and third coordinates $t$ and $\xi$ of atoms of any $N_{i,j}$, ordered according to the values of $t$. More formally, $(T_k)_{k\geqslant 1}$ is the sequence of jump times of the counting process $$\left(Z_t\right)_{t\geqslant 0}=\left(\sum_{(i,j)\in \mathcal{I}}\int_{[0,t]\times \R_+}\mathrm{1}_{z<\theta_{i,j}\left(X(u^-)\right)}N_{i,j}(\mathrm{d}u,\mathrm{d}z,[0,1])\right)_{t\geqslant 0},$$ 
so $T_k$ it is the time of the $k$-th (potentially null) transfer between patches of $\ie{1}{n}$ commanded by the $N_{i,j}$ processes. The jump quantiles sequence $(U_k)_{k\in \Z_+}$ is defined as follows:
$$\forall k\geqslant 1, \quad U_k=\int_0^1\xi \sum_{(i,j)\in \mathcal{I}}N_{i,j}(T_k,\R_+,\mathrm{d}\xi),$$
so the $U_k$ are independent and uniformly distributed on $[0,1]$, and the amplitude of the transfer from patch $i$ to patch $j$ at time $T_k$ is given by $\sum_{(i,j)\in \mathcal{I}} N_{i,j}(T_k,\R_+,U_k)q_{i,j}(X(T_k^-),U_k)$.\pe

It stems from equation (26.15) of \cite{Dav93} that the infinitesimal generator $\frak{A}$ of $\left(X,\P_x,x\in\R_+^n\right)$ is given by
\begin{equation}\label{expression_gen}\frak{A}f(x)=\sum_{k=1}^n \frac{\partial f}{\partial x_i}(x)\phi_i(x)+\sum_{(i,j)\in \mathcal{I}}\theta_{i,j}(x)\int  \bigg(f\left(x+y\cdot(e_j-e_i)\right)-f(x)\bigg)\mu_{i,j}(x,\mathrm{d}y)
\end{equation}
for every $\mathcal{C}^1$ function $f:\R_+^n\to \R$. This equation is of well-known interest in studying invariant measures for $X$, as we illustrate in Appendix A.

\section{Stability of multiplicative and unitary models}
\sectionmark{Stability results}
\label{section_models_multuni}

\indent \indent We now turn to two particular settings for which we will state simple stability results, each of them being representative of a typical feature that our graph-based model may exhibit. The first one corresponds to constant jump rates and multiplicative transfers (as in Equation (\ref{modele_motiv_multipl})) while the second one, symmetrically, is defined by coercive jump rates and unitary transfers (as in Equation (\ref{modele_motiv_unit})). Let us begin by stating additional assumptions that are common to both frameworks.

\subsection{Additional assumptions}\label{hypo_additionnelles}

\paragraph{Assumptions on autonomous growth functions}
~\ppe\\
Let $m\in \ie{1}{n}$ and $d\in \ie{1}{m}$. We set $V^+=\ie{1}{d}$, $V^0=\ie{d+1}{m}$ and $V^-=\ie{m+1}{n}$.  Patches in $V^+$, $V^0$ and $V^-$ will respectively be called \emph{sources}, \emph{neutral patches} and \emph{sinks} following the terminology of \cite{Pul}.\pe 

We consider the following set of assumptions:

\begin{hypo}[On autonomous growth functions]\label{addass_auto}
~\begin{enumerate}[label=(\alph*)]
\item For all $i\in \ie{1}{n}$ and all $x\in \R_+^n$, $\phi_i(x)=\phi_i(x_i)$ only depends on $x_i$.
\item For all $i\in \ie{1}{n}$: \begin{itemize}\item $\phi_i>0$ if $i\in V^+$ 
\item $\phi_i=0$ if $i\in V^0$ 
\item $\phi_i\leqslant 0$, $\phi_i(0)=0$ and $\phi_i(y)<0$ as soon as $y>0$ if $i\in V^-$\end{itemize}
\item For all $i\in V^+$, $\phi_i$ is continuous and piecewise $\mathcal{C}^1$.
\item For all compact subset $K\subset \R_+^n$ there exists $T>0$ such that $\Phi_i(x,T)=0$ for any $x\in K$ if $i\in V^-$.
\end{enumerate}
\end{hypo}

This first set of assumptions makes it possible to monitor the response of a given trajectory to small variations of jump times and transfer quantiles, with minor concern for the local behavior of $\Phi$, and to describe simple trajectories leading to the emptying of the system. Note that the constant growth setting, defined by 
\begin{equation}\label{constant_growth}
\phi_i(x)=\begin{cases}c_i \text{ if }c_i\geqslant 0 \\ c_i\mathrm{1}_{x_i>0} \text{ if }c_i<0\end{cases}
\end{equation} 
and involved in the multiplicative uniform and the unitary power setting, obviously satisfies these conditions.\pe 

We then require a strong condition on the graph structure of the active graph $\mathcal{G}$ to hold:
\begin{hypo}[On the topology of $\mathcal{G}$]
~\begin{enumerate}[label=(\alph*)]\label{addass_topo}
\item Any $j\in V^-$ can be reached from any $i\in \ie{1}{n}$ by a path in $\mathcal{G}$.
\item Any $j\in V^0$ can be reached from some $i\in V^+$ by a path in $\mathcal{G}$.
\end{enumerate}
\end{hypo}

\begin{figure}
\begin{center}
\begin{tikzpicture}[scale=.85,transform shape]
\tikzstyle{VertexStyle}=[shape        = circle,
                             fill         = YellowOrange,
                             minimum size = 33pt,
                             text         = black,
                             draw]
  \Vertex[L=$1$,x=-2,y=0]{A}
\tikzstyle{VertexStyle}=[shape        = circle,
                             fill         = lime,
                             minimum size = 33pt,
                             text         = black,
                             draw]
  \Vertex[L=$2$,x=-4,y=3]{B}
  \tikzstyle{VertexStyle}=[shape        = circle,
                             fill         = lime,
                             minimum size = 33pt,
                             text         = black,
                             draw]
  \Vertex[L=$3$,x=2,y=-2]{C}
  \tikzstyle{VertexStyle}=[shape        = circle,
                             fill         = lime,
                             minimum size = 33pt,
                             text         = black,
                             draw]
  \Vertex[L=$4$,x=1,y=-3]{D}
  \tikzstyle{VertexStyle}=[shape        = circle,
                             fill         = red,
                             minimum size = 33pt,
                             text         = black,
                             draw]
  \Vertex[L=$5$,x=4,y=3]{E}
  \tikzstyle{VertexStyle}=[shape        = circle,
                             fill         = red,
                             minimum size = 33pt,
                             text         = black,
                             draw]
  \Vertex[L=$6$,x=7,y=0]{F}
  \tikzstyle{EdgeStyle}=[bend right, post]
    \Edge[label=$\theta_{2,1}$](B)(A)
        \Edge[label=$\theta_{3,1}$](C)(A)
            \Edge[label=$\theta_{1,2}$](A)(B)
                         \Edge[label=$\theta_{6,5}$](F)(E)
  \tikzstyle{EdgeStyle}=[bend left, post]
        \Edge[label=$\theta_{5,3}$](E)(C)
    \Edge[label=$\theta_{1,5}$](A)(E)
        \Edge[label=$\theta_{1,6}$](A)(F)
                 \Edge[label=$\theta_{4,1}$](D)(A)
\end{tikzpicture}
\captionsetup{justification=justified}
\caption{\label{fig_ex_graphe_admissible}Graph structure complying with Assumption \ref{addass_topo}. The edges in the figure are the elements of $\mathcal{E}$; sources, neutral patches and sinks are respectively green, orange and red. Note that the graph represented is not strongly connected since no edge leads to node $4$.}
\end{center}
\end{figure}
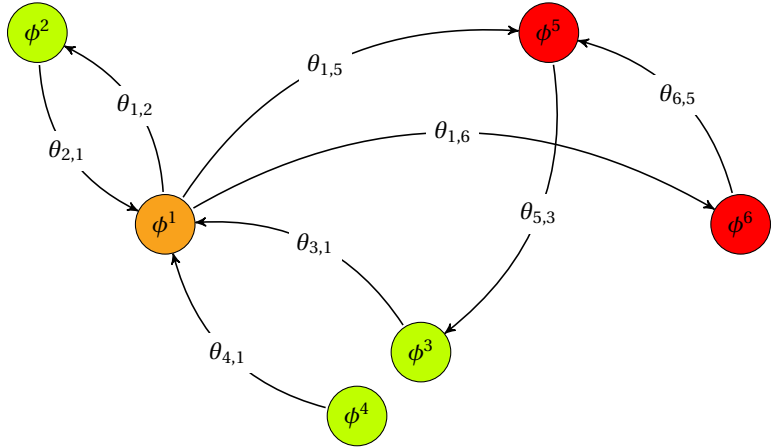

The latter assumption is quite restrictive but will prove crucial in discussing both the boundedness of the process and the petiteness of compact subsets of $\R_+^n$. It allows to describe paths along which the system empties without having to discuss complex connectivity properties, while still being weaker than plain connectivity (see Figure \ref{fig_ex_graphe_admissible}).\vspace{1.2cm}

\subsection{Multiplicative models}

We first introduce multiplicative models.

\begin{defin}[Multiplicative setting] The model is said to be \emph{multiplicative} if the following hold:
\begin{enumerate} 
\item Assumptions \ref{addass_auto} and \ref{addass_topo} are fulfilled.
\item For all $(i,j)\in \mathcal{A}$, $\theta_{i,j}$ is constant.
\item For all $(i,j)\in \mathcal{I}$, there exists a probability distribution $\mu_{i,j}$ on $[0,1]$ whose restriction to some non-punctual interval $I_{i,j}\subset [0,1]$ admits a piecewise continuous and positive density with respect to the Lebesgue measure, such that:
$$\forall x\in \R_+^n,\forall A\in \mathcal{B}([0,1]),\quad \mu_{i,j}(x,x_iA)=\mu_{i,j}\left(A\right)$$
\end{enumerate}
\end{defin}

\indent \indent The "multiplicative" denomination refers to the fact that the relative population moving from patch $i$ to patch $j$ is drawn according to a law $\mu_{i,j}$ which is independent from $x$. The process therefore goes through macroscopic jumps, but the assumption on the $\theta_{i,j}$ implies that these cannot be too frequent, which results in long, uninterrupted emptying periods under suitable assumptions. Note that we required the $\theta_{i,j}$ to be constant for the sake of simplicity, but all of the results we will present hold if the $\theta_{i,j}$ take values in some compact interval of $\R_+^*$.
\pe \\
\indent Multiplicative models encompass the introductory setting defined by Equation (\ref{modele_motiv_multipl}) as well as a large class of \emph{additive increase multiplicative decrease} (AIMD) models that are continuous counterparts to those studied in \cite{Dum} and \cite{Hon}.
\pe \\
\indent Our main stability result for the multiplicative setting is the following:
\begin{theorem}[Stability of multiplicative models]\label{thm_pcpal_mult}
If the model is multiplicative and \begin{equation}\label{cdt_flux}\limsup_{\min_i x_i \to +\infty}\sum_{i=1}^n\phi_i(x)<0\end{equation}
then $X$ is positive Harris recurrent and there exists $\eta>0$ such that $X$ is $F$-ergodic with
$$F:\begin{cases} \R_+^n\to [1,+\infty[ \\ (x_1,\ldots,x_n)\mapsto e^{\eta{\sqrt{\sum_{i=1}^n x_i}}}\end{cases}.$$
\end{theorem}

\begin{exemple}\label{ex_mod_lin}
As stated in the Introduction, a model defined by (\ref{modele_motiv_multipl}) is Harris-recurrent positive if and only if $\sum_{i=1}^n c_i<0$ and transient otherwise (remind that we assumed that $V^+\neq \varnothing$). 
\end{exemple}

\begin{figure}[h]
\begin{center}
\begin{tikzpicture}[scale=.85,transform shape]
\tikzstyle{VertexStyle}=[shape        = circle,
                             fill         = lime,
                             minimum size = 22pt,
                             text         = black,
                             draw]
  \Vertex[L=$+3$,x=-2,y=0]{A}
\tikzstyle{VertexStyle}=[shape        = circle,
                             fill         = lime,
                             minimum size = 22pt,
                             text         = black,
                             draw]
  \Vertex[L=$+2$,x=2,y=0.5]{B}
  \tikzstyle{VertexStyle}=[shape        = circle,
                             fill         = red,
                             minimum size = 22pt,
                             text         = black,
                             draw]
  \Vertex[L=$-5$,x=.5,y=-3]{C}
  \tikzstyle{VertexStyle}=[shape        = circle,
                             fill         = red,
                             minimum size = 22pt,
                             text         = black,
                             draw]
  \Vertex[L=$-1$,x=4.5,y=-2.5]{D}
  \tikzstyle{EdgeStyle}=[bend right, post]
    \Edge[label=$\theta_{1,3}$](A)(C)
        \Edge[label=$\theta_{4,3}$](D)(C)
        \Edge[label=$\theta_{3,4}$](C)(D)
  \tikzstyle{EdgeStyle}=[bend left, post]
         \Edge[label=$\theta_{1,2}$](A)(B)
                 \Edge[label=$\theta_{2,4}$](B)(D)
                 
                 \tikzstyle{VertexStyle}=[shape        = circle,
                             fill         = lime,
                             minimum size = 22pt,
                             text         = black,
                             draw]
  \Vertex[L=$+3$,x=6,y=0]{E}
\tikzstyle{VertexStyle}=[shape        = circle,
                             fill         = lime,
                             minimum size = 22pt,
                             text         = black,
                             draw]
  \Vertex[L=$+2$,x=10,y=0.5]{F}
  \tikzstyle{VertexStyle}=[shape        = circle,
                             fill         = red,
                             minimum size = 22pt,
                             text         = black,
                             draw]
  \Vertex[L=$-5$,x=8.5,y=-3]{G}
  \tikzstyle{VertexStyle}=[shape        = circle,
                             fill         = red,
                             minimum size = 22pt,
                             text         = black,
                             draw]
  \Vertex[L=$-1$,x=12.5,y=-2.5]{H}
  \tikzstyle{EdgeStyle}=[bend right, post]
   
    \Edge[label=$\theta_{1,3}$](E)(G)
        \Edge[label=$\theta_{3,4}$](G)(H)
  \tikzstyle{EdgeStyle}=[bend left, post]
             \Edge[label=$\theta_{1,2}$](E)(F)
                 \Edge[label=$\theta_{2,4}$](F)(H)
\end{tikzpicture}
\captionsetup{justification=justified}
\caption{\label{importance_struc_graphe}The importance of connectivity: if (\ref{modele_motiv_multipl}) is fulfilled, then the process corresponding to the graph configuration on the left is ergodic (according to Theorem \ref{thm_pcpal_mult}) and that on the right is transient (consider the population on the subgraph formed by the two right patches). Coefficients within patches are the $c_i$ and all $\theta_{i,j}$ represented here are strictly positive.}
\end{center}
\end{figure}
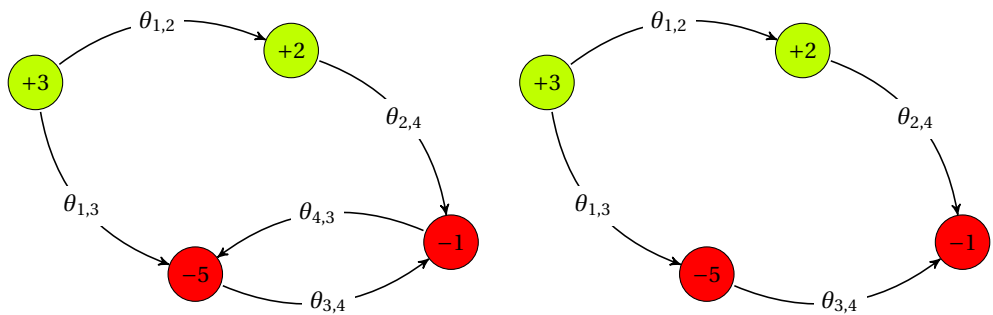

\begin{exemple}
The model represented on the left in Figure \ref{importance_struc_graphe} corresponds to an ergodic setting, while the model on the right does not.
\end{exemple}

\indent Let us give the main thrust of the proof of Theorem \ref{thm_pcpal_mult}, the details of which will be developed in a more general framework in Section \ref{section_general_crit}. Our proof strategy is based on the conceptual framework by Meyn and Tweedie \cite{MT2,MT3,DMT95}. It consists in showing that $X$ is brought back quickly enough to some compact subset of $\R_+^n$ (a \emph{boundedness} property) and that its various trajectories scan some Borel subset of $\R_+^n$ (a so-called \emph{petiteness} property). The condition on the $\phi_i$ entails the existence of some $\R_+^n$ area (namely $\{x\in \R_+^n: \min_i x_i\geqslant R\}$ for $R$ large enough) on which the flow commands a steady decrease of the system total population. The assumptions on the $\theta_{i,j}$ and the $\mu_{i,j}(x,\cdot)$ imply that with lower-bounded probability, the process reaches this zone quickly enough and stays in it long enough for the total population of the system to be brought back below a given threshold after some time, which yields the boundedness property. On the other hand, small variations of jump quantiles along a given path induce locally one-to-one changes of its point of arrival, which implies the petiteness property.\\


\indent In the transient case, the stability of transfer laws and the existence of macroscopic jumps make scaled multiplicative processes easy to describe. If the model is multiplicative with all $I_{i,j}=[0,1]$ and if $X$ is transient, then $\Vert X\Vert_1$ goes to infinity almost surely and for any $x\in \mathcal{S}=\{x\in \R_+^n\mid \Vert x \Vert_1=1\}$:
\begin{equation}\label{rel_approx_transient}\E_{Rx}\left[\sup_{t\in [0,R]}\left\Vert\frac{X(t)}{\Vert X(t)\Vert_1}-S(t)\right\Vert_{\infty}\right]\underset{R\to +\infty}{\longrightarrow} 0\end{equation}
where $S$ is a $\mathcal{S}$-valued pure jump process which is a weak solution of the following SDE:
\begin{equation*}\forall (i,j)\in \mathcal{I}, \quad \mathrm{d}S_i(t)=\sum_{j\neq i}\int_0^1 \xi S_j(t) N'_{j,i}(\mathrm{d}t,\mathrm{d}\xi)-\sum_{j\neq i}\int_0^1 \xi S_i(t) N'_{i,j}(\mathrm{d}t,\mathrm{d}\xi)\end{equation*}
with initial value $x$ under the $\P_{Rx}$ probability measure, and the $N'_{i,j}$ are independent Poisson point processes with intensity $\mathrm{d}t\mu_{i,j}(\mathrm{d}\xi)$. Studying the stability of $S$ only requires petiteness analysis, and (\ref{rel_approx_transient}) makes it possible to infer the behavior of $\frac{X}{\Vert X \Vert_1}$ (which may not be a Markov process) from that of $S$. As an example, we can prove the following result (see Appendix B):

\begin{framed}\vspace*{-.45cm}\begin{prop}[Transient multiplicative uniform case when $n=2$]\label{gamma_transient}
In the two-patch multiplicative uniform model with $c_1+c_2\geqslant 0$, $\frac{X_1(t)}{X_1(t)+X_2(t)}$ converges in distribution to a Beta distribution with parameters $\left(\frac{\theta_{2,1}}{\lambda},\frac{\theta_{1,2}}{\lambda}\right)$ as $t$ goes to infinity, where $\lambda=\theta_{1,2}+\theta_{2,1}$.
\end{prop}\vspace*{-.45cm}\end{framed}

\subsection{Unitary models}

We define unitary models as follows:

\begin{defin}[Unitary setting]The PDMP metapopulation model is said to be \emph{unitary} if the following hold:
\begin{enumerate} 
\item Assumptions \ref{addass_auto} and \ref{addass_topo} are fulfilled.
\item There exists an increasing and subadditive function $\Theta:\R_+\to \R_+^*$ with $\lim_{y\to +\infty}\Theta(y)=+\infty$ and $\alpha>1$ such that for any $x\in \R_+^n$:
\begin{itemize}
\item For all $(i,j)\in \mathcal{E}$, $\Theta(x_i)\leqslant \theta_{i,j}(x)$.
\item For all $(i,j)\in \mathcal{I}$, $\theta_{i,j}(x)\leqslant \alpha\Theta(x_i)$.
\end{itemize}
\item For all $x\in \R_+^n$ and all $(i,j)\in \mathcal{E}$, $\mu_{i,j}(x,\cdot)=\delta_{1\wedge x_i}$.\vspace{1mm}~
\end{enumerate}
\end{defin}

Our main stability result on unitary models is the following:

\begin{theorem}\label{thm_pcpal_uni}
If the model is unitary and if either
\begin{equation}\label{flux_strong}\limsup_{\min_{i\in V^-} x_i \to +\infty}\sum_{i=1}^n\phi_i(x)<0\end{equation}
or \begin{equation}\label{flux_weak}{\mathcal{G}\text{ is strongly connected and}}\limsup_{\min_i x_i \to +\infty}\sum_{i=1}^n\phi_i(x)<0\end{equation}
then there exists $\eta>0$ such that $X$ is positive Harris recurrent and $F$-ergodic with
$$F:\begin{cases} \R_+^n\to [1,+\infty[ \\ (x_1,\ldots,x_n)\mapsto e^{\eta{\sqrt{\sum_{i=1}^n x_i}}}\end{cases}.$$
\end{theorem}

\begin{exemple}
The following setting corresponds to an unitary metapopulation model with logistic autonomous growth in sources and jump rates taking into account the carrying capacity of target patches:
\begin{equation}\label{toy_model_log}\begin{cases}
\forall i\in V^+,\quad \phi_i(x)=b_i x_i(K_i-x_i)_++c_i \\
\forall j\in V^-,\quad \phi_j(x)=-\frac{c_jx_j}{\alpha_j+x_j}\\
\forall i\in V^+,\forall j\in V^-, \quad \theta_{i,j}(x)=\gamma_{i,j}x_i\frac{\varepsilon_{i,j}+x_j}{1+x_j}\quad \text{and}\quad \theta_{j,i}\equiv 0\\
\forall (i,j)\notin \mathcal{E},\quad \theta_{i,j}\equiv 0\\
\forall (i,j)\in \mathcal{E},\quad \mu_{i,j}(x,\cdot)=\delta_{1\wedge x_i}
\end{cases}\end{equation}
where the $b_i$, $K_i$, $\gamma_{i,j}$ and $\varepsilon_{i,j}$ are positive constants and the $c_i$ are such that $c_i>0$ if $i\in V^+$ and $c_i<0$ if $i\in V^-$. Theorem \ref{thm_pcpal_uni} implies the ergodicity of the process in this setting as soon as $\sum_{i=1}^n c_i < 0$ if $\mathcal{G}$ is strongly connected.
\end{exemple}

\begin{exemple}
As stated in our Introduction, a power unitary model defined by (\ref{modele_motiv_unit}) is ergodic if and only if $\sum_{i=1}^n c_i<0$ and transient otherwise (remind that we assumed that $V^+\neq \varnothing$).
\end{exemple}

The proof of Theorem \ref{thm_pcpal_uni} is based on arguments that only slightly differ from those put forward in the multiplicative case. If $\mathcal{G}$ is strongly connected, boundedness is shown as in the proof of Theorem \ref{thm_pcpal_mult}. If this is not the case but (\ref{flux_strong}) holds, we show that $X$ satisfies Assumption \ref{majo_temps_retour_zone_decr} for any value of $T$, which entails that Assumption \ref{effic_phases_descente} holds as well. In both cases, we then consider small variations of jump times (and not jump quantiles anymore) to check for the petiteness of compacts for the resolvent chain of $X$.

\section{Stability criteria for the general metapopulation model}
\label{section_general_crit}
\indent \indent We now present Meyn-Tweedie inspired criteria for boundedness, petiteness and ergodicity for the general model defined in Section \ref{section_modelisation}. We will see that their application requires proof strategies that are dependent from the model specification and the system's active graph structure. However, it is easy to adapt the proofs of Theorems \ref{thm_pcpal_mult} and \ref{thm_pcpal_uni} in a discretionary way to a large range of frameworks, which is why we endeavoured to state their general philosophy. In particular, these criteria will apply to both multiplicative and unitary settings, yielding Theorems \ref{thm_pcpal_mult} and \ref{thm_pcpal_uni}.\pe 

\subsection{Criterion for boundedness}
\label{section_bornitude}

\indent \indent We first want to state a sufficient condition for $X$ to be \emph{bounded in probability on average}. Remind from \cite{MT2} that this property means that for all $x\in \R_+^n$ and all $\varepsilon >0$ there exists a compact subset $C\subset \R_+^n$ such that $$\liminf_{t\to +\infty} \frac{1}{t}\int_0^t \P_x(X(s)\in C)\mathrm{d}s\geqslant 1-\varepsilon.$$

Fix a Borel subset $S$ of $\R_+^n$. Our first assumption implies that the flow on $S$ drives the process at a steady rate towards the origin:

\begin{hypo}[Existence of a steady population decrease zone]\label{hyp_zone_decr}
There exists $c>0$ such that:
$$\forall x\in S, \quad \sum_{i=1}^n \phi^i(x)\leqslant -c$$
\end{hypo}

\noindent $S$ corresponds to a set of population configurations in which the local deterministic dynamics cause the system to empty. For instance, in the constant growth setting defined by (\ref{constant_growth}) and under condition (\ref{cond_flux_tot}), $S$ may be defined as any Borel set of population configurations in which all sinks contain positive population.\pe

\noindent It is worth mentioning a simple case with straightforward consequences. If Assumption \ref{hyp_zone_decr} holds and with $S=K^c$ for some compact subset $K$ of $\R_+^n$, it is a simple matter to show that $X$ is bounded in probability on average using, for instance, Fatou's Lemma. In this case, $\sum \phi_i$ may play the role of a Lyapunov function for $X$, and it is possible to derive results on recurrence and ergodicity provided all compacts are petite and there exists an irreducible skeleton chain for $X$ (see Theorem 4.2 of \cite{MT3} and Theorem 5.2 of \cite{DMT95}). However, it is not always the case that $S=K^c$ for some compact $K$. In particular, it is not under Assumption \ref{addass_auto} since $S$ cannot contain any $x$ such that $x_i=0$ for all $i\in V^-$.\pe

\noindent Under Assumption \ref{hyp_zone_decr}, it is natural to think that $X$ will be bounded in probability on average if it quickly goes back to $S$ and stays within $S$ for a long time whenever it reaches it. This is why we put forward the following conditions:

\begin{hypo}[Bounds for the hitting time of $S'$ and the exit time from $S$]\label{assump_bound}
\pe \\ There exists a Borel subset $S'$ of $S$ as well as $\delta,\varepsilon,T,T'\in \R_+^*$ and $R\geqslant 0$ such that:
\begin{enumerate}
\item \label{majo_temps_retour_zone_decr} For any $x\notin S$ with $\Vert x \Vert_1\geqslant R$:
\begin{equation*}\P_x(\exists s\in [0,T], X_s\in S')\geqslant \delta \end{equation*}
\item 
 \label{mino_temps_depart_zone_decr}For any $x\in S'$ with $\Vert x \Vert_1\geqslant R$: 
\begin{equation*}\P_x(\forall s\in [0,T'],X_s\in S)\geqslant \varepsilon\end{equation*}
\item \label{effic_phases_descente} We have \begin{equation*}\varepsilon T' c > (1-\varepsilon)\frac{T}{\delta}\sup_{\R_+^n}\sum_{i=1}^n \phi^i\end{equation*}
\end{enumerate}
\end{hypo}

\indent Note that Assumption \ref{assump_bound}.\ref{mino_temps_depart_zone_decr} may be void since we do not assume that $S'\cap \{\Vert x \Vert_1\geqslant R\}$ is not empty.\pe

\indent One can see that under Assumptions \ref{hyp_zone_decr}, \ref{assump_bound}.\ref{majo_temps_retour_zone_decr} and \ref{assump_bound}.\ref{mino_temps_depart_zone_decr} and starting from a point in $S$ with high enough total population, then with probability of at least $\varepsilon$, $X$ does not leave $S$ before time $T'$, and the mean total population increase during an excursion of the process outside of $S$ is at most $\frac{T}{\delta}\sup_{\,\R_+^n}\sum_i\phi_i$. Assumption is \ref{assump_bound}.\ref{effic_phases_descente} therefore sufficient for the visits of $S$ by $X$ to bring the population process back to some compact set of $\R_+^n$ regardless of its original position.
\pe \\
The expected result follows. Its proof is deferred to Section \ref{section_demos} and consists in a comparison with a random walk on $\R$. It may easily adapt to some models that do not necessarily fulfill Assumptions \ref{hyp_zone_decr} and \ref{assump_bound} by considering suitable $\R$-valued Markov chains.

\begin{theorem}\label{thm_born_prob}
Under Assumptions \ref{hyp_zone_decr} and \ref{assump_bound}, $X$ is bounded in probability on average.
\end{theorem}

Alternate criteria for boundedness may be derived using Theorem 2.1 of \cite{MT4} under irreducibility assumptions.  In particular, Theorem 3.1 of \cite{Dai} holds in our setting whenever compact subsets of $\R_+^n$ are petite. Yet, as we will see in Section \ref{section_ergo}, Assumptions \ref{hyp_zone_decr} and \ref{assump_bound} yield an upper bound for expected exponential functionals of return times, which implies strong ergodicity results.

\subsection{Criterion for petiteness}
\label{section_petitesse}

\indent \indent Let us now consider a compact set $C\subset \R_+^n$. We remind (see \cite{MT2}) that $C$ being \emph{petite} for some sampled chain $(X(A_n))_{n\geqslant 0}$ of $X$ means that there exists a non-trivial Borel measure $\nu$ on $\R_+^n$ such that $\P_x\left(X(A_1)\in \cdot \right)\geqslant \nu(\cdot)$ for all $x\in C$, and that $C$ is just said to be \emph{petite} if it is petite for some $(A_n)_{n\geqslant 0}$.
\pe \\
\indent Proving that $C$ is petite relies on framework-specific strategies. The criterion for petiteness we will derive in this section applies to a broad range of settings, but it is worth keeping in mind that it may prove unnecessarily technical in some cases. In particular, it can be easier to identify sampled chains that dominate $\R_+^n$-valued Dirac measures --- provided, of course, that such chains exists.\pe 

In most non-pathological specifications, the distribution of inter-jump times exhibits an absolutely continuous component with respect to the Lebesgue measure. This observation drives us to look for subsets of $\R_+^n$ on which the Lebesgue measure is dominated by the semigroup for the resolvent chain $(R_n)_{n\in \Z_+}$ defined by:
$$\forall n\in \Z_+,\quad R_n=X(S_n)$$
where $(S_n)_{n\in \Z_+}$ is the sequence of jump times of a Poisson process with density $1$ independent from $X$.\pe 

With this aim in mind, let us first define tracking functions that return the position of our process from its inter-jump times $t_k$, the edges $\kappa_k$ along which the transfers occur and the values of its jump distributions quantiles $\xi_k$. For all $x\in \R_+^n$, all $(i,j)\in \mathcal{I}$ and all $\theta\in [0,1]$, set
$$g_{i,j}(x,\xi)=x+q_{i,j}\left(x,\xi\right)(e_j-e_i).$$
For all $x\in \R_+^n$, let us define recursively the sequence of functions $\left(h^k_x\right)_{k\geqslant 1}$ by
$$h_x^1:\begin{cases} \bigcup_{l\geqslant 1} \left(\R_+^l\times [0,1]^{l-1}\times \mathcal{I}^{l-1}\right)\to \R^n \\ (t,\xi,\kappa)\mapsto \Phi(x,t_1)\end{cases}$$ 
and for all $k\in \N$, 
$$h_x^{k+1}:\begin{cases} \bigcup_{l\geqslant k+1} \left(\R_+^l\times [0,1]^{l-1}\times \mathcal{I}^{l-1}\right) \to \R^n \\ (t,\xi,\kappa) \mapsto \Phi\left(g_{\kappa_k}\left(h_x^k(t,\xi,\kappa),\xi_k\right),t_{k+1}\right)\end{cases}$$
where $\kappa=(\kappa_1,\ldots,\kappa_{j-1})$ denotes the generic element of $\mathcal{I}^{j-1}$. The vector $h_x^k(t,\xi,\kappa)$ is the state occupied by the process initiated at $x$ after it has followed the flow for time $t_1$, undergone a transfer along edge $\kappa_1$ with amplitude given by the quantile of order $\xi_1$ of the appropriate $\mu_{i,j}(x,\cdot)$ law, followed the flow again for time $t_2$, and so on until $k-1$ transfers occurred and the process followed the flow for time $t_k$ after its last jump.\pe

Our strategy is to look for a subset $V$ of an affine subspace of $\R_+^n$ such that we can provide a lower bound for the Lebesgue measure of the pre-image by some $h_x^k$ of any Borel subset of $V$. We therefore introduce the following assumption:

\begin{hypo}[Likely paths scanning a Borel subset]\label{assump_petite}\pe \\
There are $M>0$, $\overline{N}\in \N$, $\overline{T}\geqslant 1$, $p\in \ie{1}{n}$, a $p$-dimensional affine subspace $V$ of ~$\R^n$, a Borel subset $\mathcal{P}$ of\, $V$ with non-zero Lebesgue measure in $V$, and for any $x\in C$ there are $N(x)\in \ie{1}{\overline{N}}$ with $q(x)=2N(x)-1-p\geqslant 0$, a vector of edges $\kappa_x\in \mathcal{I}^{N(x)-1}$, a\, $\R^{2N(x)-1}$-coordinates permutation $\sigma_x$ and open subsets $U_1^x\subset [0,\overline{T}]^p$ and $U_2^x\subset [0,\overline{T}]^{q(x)}$ such that setting
$$\psi_{z_2}^x:\begin{cases} U_1^x\to V \\ z_1\mapsto h_x^{N(x)}(\sigma_x(z_1,z_2),\kappa_x)\end{cases} $$
the following hold:
\begin{enumerate} 
\item \label{hyp_dilat} For all $z_2\in U_2^x$, $\psi^x_{z_2}$ is a $\mathcal{C}^1$-diffeomorphism with Jacobian determinant bounded by $M$ on $U_1^x$.
\item \label{constant_scanned} For all $z_2\in U_2^x$, $$\mathcal{P}\subset \psi^x_{z_2}(U_1^x).$$
\item 
\label{hyp_mino_dens_sauts} Defining $$\theta_0=\inf_{x\in C}\left[\lambda_{q(x)}(U_2^x)\inf_{\substack{z_1\in U_1^x \\z_2\in U_2^x}}\prod_{i=1}^{N(x)-1}\theta_{\kappa_i(x)}\left(h_x^k(\sigma_x(z_1,z_2),\kappa_x)\right)\right],$$
we have that $\theta_0>0$.
\end{enumerate}
\end{hypo}

\noindent Note that Proposition \ref{prop_petits} below still holds if we replace $p$-dimensional affine subspaces by $p$-dimensional manifolds in the condition above. Condition \ref{assump_petite}.\ref{hyp_mino_dens_sauts} states that paths that lead $X$ to Borel subsets of $V$ do not correspond to unlikely sequences of jump times or jump quantiles.
\pe \\
\indent A change of variables argument then yields the following proposition:
\begin{propgrise}\label{prop_petits}
If Assumption \ref{assump_petite} is met, then $C$ is petite.
\end{propgrise}

Most of the technicity in applying this proposition lies on proving that Assumption \ref{assump_petite}.\ref{hyp_dilat} holds. Indeed, this requires describing paths of the process that lead to a given area of the state space as well as monitoring their response to small variations of jump times and jump quantiles. In general, this can be made simpler by the following straightforward but useful result allowing for localization of the starting point of such paths, which is Lemma 3.1 from \cite{MT1}.

\begin{framed}\vspace*{-.45cm}\begin{prop}[Petiteness transitivity]\label{transmission_petitesse}
If $C'$ is a Borel subset of $\R_+^n$ and if there exists $A>0$ such that
$$\inf_{x\in C}~\P_x\left(\tau_{C'}<A\right)>0,$$ then if $C$ a petite set for the resolvent of $X$, so is $C$.
\end{prop}\vspace*{-.45cm}\end{framed}

%

\subsection{Criterion for ergodicity}
\label{section_ergo}

\indent \indent Theorems 3.2 (ii) and 4.1 (i) of \cite{MT2} along with present Theorem \ref{thm_born_prob} and Proposition \ref{prop_petits} imply that $X$ is positive Harris recurrent as soon as Assumptions \ref{hyp_zone_decr} and \ref{assump_bound} are met and all compact subsets of $\R_+^n$ are petite. We now look for an additional condition to ensure \emph{$F$-ergodicity} for some measurable $F:\R_+\to [1,+\infty[$. Remind from \cite{MT2} that this property writes:
$$\forall x\in \R_+^n,\quad \lim_{t\to +\infty}\sup_{|g|\leqslant F}\left\vert \E_x\left(g(X(t))\right)-\pi(g)\right\vert=0$$ 
with $\pi$ standing for the invariant probability of $X$ and the supremum being taken over all measurable $g:\R_+^n\to \R$ functions such that $|g|\leqslant F$.
\pe \\
Our result is the following:
\begin{theorem}\label{thm_ergo}
Assume that Assumptions \ref{hyp_zone_decr} and \ref{assump_bound}, that all compact subsets of $\R_+^n$ are petite and that $X$ admits an irreducible skeleton chain. Then there exists $\eta\in \R_+^*$ such that $X$ is $F$-ergodic for 
\begin{equation}\label{defF}F:\begin{cases} \R_+^n\to \R_+\\ x\mapsto e^{\eta \sqrt{\sum_{i=1}^n x_i}}\end{cases}\end{equation}
\end{theorem}
Theorems \ref{thm_pcpal_mult} and \ref{thm_pcpal_uni} are corollaries of the latter result, as is detailed in Section \ref{section_demos}.

\section*{Conclusion}

 \indent  \indent Besides metapopulations, possible fields of application for the class of piecewise deterministic models presented here range from open Jackson networks \cite{Jack,Mjack,Dai} with deterministic inputs and outputs to communication networks based on TCP-type processes \cite{TCP,TCP2,Hon}, storage \cite{Har} on a network, neuronal stimulation \cite{Ost,Dua} and more generally a large class of stochastic hybrid systems on graphs \cite{Cas} with low stochasticity in autonomous dynamics.\ppe

\indent Our work may be expanded in many ways. First, the multiplicative and unitary frameworks were designed as simple models that allow either for large jumps or large jump rates, but many applications may require designing and studying hybrid models. Criteria from Sections \ref{section_bornitude} and \ref{section_petitesse} will hopefully prove flexible and be useful in such settings.\ppe 

\indent Besides, some metapopulation settings do not necessarily fit with the additional hypothesis we made on the structure of the active graph $\mathcal{G}$ of the system --- that is, on the graph formed by edges along which non-zero transfers occur at non-zero rate. Although the criteria for boundedness and petiteness we stated in Sections \ref{section_bornitude} and \ref{section_petitesse} do not refer to this active graph, they cannot be applied without particular sequences of transfers being made explicit, which assumes that one can describe entire paths followed by the population load. Moreover, they require some degree of connectivity so the system, loosely speaking, can "empty" and "mix". Expanding the results above to more general graphs is still a work in progress.\ppe

\indent Finally, we chose to model deterministic intra-patch population dynamics. This assumption is only legitimate if the local demographics exhibit little stochasticity, either intrinsically or because quantities involved are so large that their aggregate evolution can be proxied by a non-random dynamic system. If this condition is not met, intra-patch birth and death or Hawkes-inspired processes (\cite{Haw}, \cite{Del}) should be considered. While this would require developing new tools for petiteness analysis, we are confident of our boundedness results holding without any major alteration.

\section*{Acknowledgements}

\indent \indent I express my deepest gratitude to Vincent Bansaye and Elisabeta Vergu for their guidance during the conception of this paper. This work was supported by the French Research Agency within projects ANR-16-CE32-0007-01 (CADENCE) and ANR-16-CE40-0001 (ABIM), and by Chaire Modélisation Mathématique et Biodiversité Veolia-X-MNHM-FX.

\section{Proofs}
\label{section_demos}

\indent \indent This section contains the proofs of the results stated previously. We begin by setting up some notation.
\\
We write $$M:=\sup_{\R_+^n}\sum_{i=1}^n \left|\phi_i\right|<+\infty.$$
For any Borel subset $E\subset\R_+^n$, we denote by $\left(\tau^k_E\right)_{k\geqslant 1}$ the sequence of successive hitting times of $E$ by $X$ (with $\tau^1_E=0$ $\P_x$-a.s. if $x\in E$), so $\tau_E=\tau^1_E$.\pe 

\subsection{Boundedness criterion \ref{thm_born_prob}}

\begin{demthm}[\ref{thm_born_prob}]We proceed in three steps to prove Theorem \ref{thm_born_prob}. 
\paragraph*{Step 1: study of a random walk}
~\pe \\
\indent Assume that Assumptions \ref{hyp_zone_decr} and \ref{assump_bound} hold and let $y\in \R$. Consider a $\R$-valued random walk $Y=\left(Y_k\right)_{k\geqslant 1}$ with i.i.d. increments defined on $(\Omega,\mathcal{A},\P)$ in such a way that: \begin{equation}\label{eq_marcheal}Y_1=y \quad \text{and} \quad \forall k\geqslant 1, \quad Y_{k+1}-Y_k\overset{d}{=}-BcT'+(1-B)\Gamma \cdot TM\end{equation}
where $B$ and $\Gamma$ are independent random variables with respective distributions Bernoulli $\mathcal{B}(1,\varepsilon)$ and shifted geometric $\mathcal{G}(\delta)$. We first observe that
$$\E\left(-BcT'+(1-B)\Gamma \cdot TM\right)=-\varepsilon c T' +(1-\varepsilon)\frac{T}{\delta}M<0$$ which entails that the hitting time $\sigma_{R}$ of $(-\infty,R]$ by $Y$ is almost surely finite by the law of large numbers. Moreover, there exists $r>0$ independent from the choice of $y\in \R$ such that:
$$\E\left(e^{\gamma(r)\sigma_{R}}\right)\leqslant e^{r\left(-cT'-R+y\right)}.$$
Indeed, straightforward computations show that if $r>0$ is such that $r<\frac{\left|\log(1-\delta)\right|}{TM}$, and if we set $$\gamma(r)=-\log\left(\frac{\delta(1-\varepsilon)}{e^{-rTM}-(1-\delta)}+\varepsilon e^{-rc T'}\right)$$ then $\left(\exp\left(rY_k+\gamma(r)(k-1)\right)\right)_{k\geqslant 1}$ is a positive martingale with respect to its natural filtration. Using Fatou's lemma then yields the expected inequality since $\gamma(r)>0$.

\paragraph*{Step 2: exponential moment of the hitting time of a compact subset}
~\pe \\
\indent Let $C=\{x\in \R_+^n\mid \Vert x \Vert_1\leqslant R\}$. We now prove that there exists $\beta>\R_+^*$ such that 
$x\mapsto \E_x\left(\exp\left(\beta\sqrt{\tau_C}\right)\right)$ is locally bounded on $\R_+^n$.\pe  

Let us first assume that $x\in S^c\cap C^c$. Then $\P_x$-a.s.:
\begin{align*}\tau_C ~\leqslant ~& \sum_{k=1}^{+\infty}\mathrm{1}_{\tau^k_{S^c}<\tau_C}\left(\tau^{k+1}_{S^c}-\tau^k_{S^c}\right) = \sum_{k=1}^{+\infty}\mathrm{1}_{\tau^k_{S^c}<\tau_C}\left(\tau^{k}_{S}-\tau^k_{S^c}+\tau^{k+1}_{S^c}-\tau^{k}_{S}\right)&\\ ~\leqslant~& \sum_{k=1}^{+\infty}\mathrm{1}_{\tau^k_{S^c}<\tau_C}\left(\tau_S^{k}-\tau_{S^c}^k+\frac{\Vert X(\tau_{S}^{k})\Vert_1-R}{|c|}\right)&\end{align*}
since the decrease of the flow on $S\cap C^c$ is at least $|c|$. For $\beta>0$, we thus get:
\begin{multline}\label{eq_demo_temps_retour_int}\left(\E_x\left[\exp\left(\beta \sqrt{\tau_C}\right)\right]\right)^2 \\ \leqslant \E_x\left[\exp\left(2\beta \sqrt{\sum_{k=1}^{+\infty}\mathrm{1}_{\tau^k_{S^c}<\tau_C}\left(\tau_S^{k}-\tau_{S^c}^k\right)} \right)\right]\E_x\left[\exp\left(2\beta \sqrt{\sum_{k=1}^{+\infty}\mathrm{1}_{\tau^k_{S^c}<\tau_C}\frac{\Vert X(\tau_{S}^{k})\Vert_1-R}{c}} \right)\right]\end{multline}
$\P_x$-a.s. using the inequality $\sqrt{a+b}\leqslant \sqrt{a}+\sqrt{b}$ for $(a,b)\in \R_+^2$ and the Cauchy-Schwarz inequality. Yet Assumption \ref{assump_bound} entails that the increment of the total system population each time the process leaves $S$ before $\tau_C$ is dominated by $-BcT'+(1-B)\Gamma \cdot TM$. Applying the strong Markov property to the sequences of stopping times $\left(\tau^k_S\wedge \tau_C\right)_{k\geq 1}$ and $\left(\left(\tau^k_{S^c}+pT\right)\wedge \tau^{k+1}_{S}\wedge \tau_C\right)_{k\geq 1,p\geq 0}$ ensures that there exists $\left(Y_k\right)_{k\geqslant 1}$ on $(\Omega,\mathcal{A},\P_x)$ satisfying (\ref{eq_marcheal}) with $y=\Vert x \Vert_1$ such that $\P_x$-a.s.:
\begin{equation}\label{eq1_domsto} \forall k\geqslant 1,\quad \mathrm{1}_{\tau^k_{S^c}<\tau_C}\left(\tau_S^k-\tau_{S^c}^k\right)\leqslant \mathrm{1}_{k<\sigma_{R}}\left(\frac{Y_{k+1}-Y_k}{M}+cT'\right)\end{equation}
and
\begin{equation}\label{eq2_domsto}\forall k\geqslant 1,\quad \left\Vert X(\tau^k_{S^c}\wedge \tau_C)\right\Vert_1 \leqslant Y_{k\wedge \sigma_{R}}.\end{equation}
For any $\beta>0$, (\ref{eq1_domsto}) yields
$$\E_x\left[\exp\left(2\beta \sqrt{\sum_{k=1}^{+\infty}\mathrm{1}_{\tau^k_{S^c}<\tau_C}\left(\tau_S^{k}-\tau_{S^c}^k\right)} \right)\right] \leqslant \E\left[\exp\left(2\beta\sqrt{\sum_{k=1}^{\sigma_{R}-1}\left(\frac{Y_{k+1}-Y_k}{M}+cT'\right)}\right)\right]$$
so \begin{equation}\label{majo_esp_1}\E_x\left[\exp\left(2\beta \sqrt{\sum_{k=1}^{+\infty}\mathrm{1}_{\tau^k_{S^c}<\tau_C}\left(\tau_S^{k}-\tau_{S^c}^k\right)} \right)\right]\leqslant \E\left[\exp\left(2\beta\sqrt{|c|T'\sigma_{R}}\right)\right].\end{equation}
Now (\ref{eq2_domsto}) entails:
\begin{align*}\E_x\left[\exp\left(2\beta \sqrt{\sum_{k=1}^{+\infty}\mathrm{1}_{\tau^k_{S^c}<\tau_C}\frac{\Vert X(\tau_{S}^{k})\Vert_1-R}{c}} \right)\right] ~\leqslant~  \E\left[\exp\left(2\beta\sqrt{\sum_{k=1}^{\sigma_{R}-1}\frac{Y_k}{c}}\right)\right].\end{align*}
The increments of the $Y$ chain are greater than $-ct$ and the value of $Y$ at time $\sigma_{R}$ is at least $R$, so we can write
$$\sum_{k=1}^{\sigma_{R}-1}Y_{k}\leqslant \sum_{k=1}^{\sigma_{R}-1}\left(R+\sigma_{R}cT'\right)\leqslant\sigma_{R}R+\sigma_{R}^2cT'$$
from which we deduce that
\begin{equation}\label{majo_esp_2}\E_x\left[\exp\left(2\beta \sqrt{\sum_{k=0}^{+\infty}\mathrm{1}_{\tau^k_{S^c}<\tau_C}\frac{\Vert X(\tau_{S}^{k})\Vert_1-R}{c}} \right)\right]\leqslant \E\left[\exp\left(2\beta\sqrt{\frac{R}{c}\sigma_{R}+T'\sigma_{R}^2}\right)\right].\end{equation}
Inequality (\ref{eq_demo_temps_retour_int}) combined with (\ref{majo_esp_1}) and (\ref{majo_esp_2}) finally yields:
\begin{equation}\label{eq_2_demo_temps_retour_int}
\left(\E_x\left(e^{\beta \sqrt{\tau_C}}\right)\right)^2\leqslant \E\left[\exp\left(2\beta\sqrt{cT'\sigma_{R}}\right)\right]\E\left[\exp\left(2\beta\sqrt{\frac{R}{c}\sigma_{R}+T'\sigma_{R}^2}\right)\right].
\end{equation}
Obviously $\sigma_{R}\leqslant \sigma_{R}^2$, so there exist $\alpha_1,\alpha_2\in \R_+^*$ independent from our choice of $x$ in $S^c\cap C^c$ such that
\begin{equation*}\left(\E_x\left(e^{\beta \sqrt{\tau_C}}\right)\right)^2\leqslant \E\left[\exp\left(2\beta\alpha_1\sigma_{R}\right)\right]\E\left[\exp\left(2\beta\alpha_2\sigma_{R}\right)\right]
\end{equation*}
and therefore
\begin{equation*}\E_x\left(e^{\beta \sqrt{\tau_C}}\right)\leqslant \E\left[\exp\left(2\beta(\alpha_1\vee \alpha_2)\sigma_{R}\right)\right].
\end{equation*}
We derive a similar inequality from
$$\tau_C ~\leqslant ~ \frac{\Vert x \Vert_1 - R}{c}+\sum_{k=1}^{+\infty}\mathrm{1}_{\tau^k_{S^c}<\tau_C}\left(\tau^{k}_{S}-\tau^k_{S^c}\right)$$
if $x\in S\cap C^c$, and the case $x\in C$ is trivial. Step 1 ensures that if $\beta>0$ is small enough, then the function $x\mapsto \E_x\left(e^{\beta\sqrt{\tau_C}}\right)$ is locally bounded on $\R_+^n$. We will use this result when proving the ergodicity of $X$ in section \ref{section_ergo}; for our present purpose, it is sufficient to know that $x\mapsto \E_x\left[\left(\tau_C\right)^2\right]$ is locally bounded on $\R_+^n$.

\paragraph*{Step 3: boundedness in probability on average}
~\pe \\
\indent If $p\geqslant 3$ is an integer, then for all $t>0$ and $x\in \R_+^n$ :
$$\frac{1}{t}\int_0^t \mathrm{1}_{X(s)\notin pC}\,\mathrm{d}s\leqslant \frac{1}{t}\left[\tau_C+\sum_{j\geqslant 1,k\geqslant 1}\mathrm{1}_{\tau_C^k<\tau_{(2C)^c}^j<\tau_C^{k+1}}\mathrm{1}_{\tau_{(2C)^c}^j<t}\int_{\tau_{(2C)^c}^j}^{\tau_C^{k+1}} \mathrm{1}_{X(s)\notin pC}\,\mathrm{d}s\right]$$
since $(pC)^c\subset (2C)^c$, so that
$$\E_x\left[\frac{1}{t}\int_0^t \mathrm{1}_{X(s)\in pC}\,\mathrm{d}s\right]\leqslant \frac{1}{t}\E_x\left[\tau_C+\sum_{j,k\in \Z_+^*}\mathrm{1}_{\tau_C^k<\tau_{(2C)^c}^j<\tau_C^{k+1}}\mathrm{1}_{\tau_{(2C)^c}^j<t}\left(\tau_C^{k+1}-\tau_{(2C)^c}^j\right)\mathrm{1}_{\tau_C^{k+1}-\tau_{(2C)^c}^j\geqslant \frac{(p-2)R}{M}}\right]$$
because the process needs at least $\frac{(p-2)R}{M}$ units of time to reach $pC$ from a state with total population $2R$. From this we deduce
$$\E_x\left[\frac{1}{t}\int_0^t \mathrm{1}_{X(s)\in pC}\,\mathrm{d}s\right]\leqslant \frac{1}{t}\E_x\left[\tau_C+\sum_{j,k\in \Z_+^*}\mathrm{1}_{\tau_C^k<\tau_{(2C)^c}^j<\tau_C^{k+1}}\mathrm{1}_{\tau_{(2C)^c}^j<t}\E_{X_{\tau_{(2C)^c}^j}}\left(\tau_C\mathrm{1}_{\tau_C\geqslant \frac{(p-2)R}{M}}\right)\right]$$
using the strong Markov property, and, setting $\zeta=\sup_{y\in \R_+^n,\Vert y \Vert_1\leqslant 2C}\E_y\left[\left(\tau_C\right)^2\right]$ (which is finite according to Step 2) and writing that $\tau_C\mathrm{1}_{\tau_C\geqslant \frac{(p-2)R}{M}}\leqslant \frac{M}{(p-2)R}\left(\tau_C\right)^2$,
$$\E_x\left[\frac{1}{t}\int_0^t \mathrm{1}_{X(s)\in pC}\,\mathrm{d}s\right]\leqslant \frac{1}{t}\left(\zeta+\E_x\left[\sum_{j,k\in \Z_+^*}\mathrm{1}_{\tau_C^k<\tau_{(2C)^c}^j<\tau_C^{k+1}}\mathrm{1}_{\tau_{(2C)^c}^j<t}\zeta \sqrt{\frac{M}{(p-2)R}}\right]\right)$$
by the Cauchy-Schwarz inequality. Thus, for any integer $p\geqslant 3$, any $t>0$ and any $x\in \R_+^n$:
$$\E_x\left[\frac{1}{t}\int_0^t \mathrm{1}_{X(s)\in pC}\,\mathrm{d}s\right]\leqslant \frac{1}{t}\left(\zeta+\frac{t}{MR}\zeta \sqrt{\frac{M}{(p-2)R}}\right)$$
as the process cannot go through more than $\frac{t}{MR}$ times the full way between $C$ and $(2C)^c$ within time $t$. Finally,
$$\liminf_{t\to +\infty}~\E_x\left[\frac{1}{t}\int_0^t \mathrm{1}_{X(s)\in pC}\mathrm{d}s\right]\leqslant \frac{\zeta}{\sqrt{M}R^{3/2}} \sqrt{\frac{1}{p-2}}.$$
Choosing $p$ arbitrarily large shows that $X$ is bounded in probability on average.\hfill $\square$\end{demthm}

\subsection{Petiteness criterion \ref{prop_petits}}

\begin{demprop}[\ref{prop_petits}]
Let us suppose that Assumption \ref{assump_petite} is met and recall that we denote the resolvent of $X$ by $(R_n)_{n\in \N}$. For all $x\in C$ and all $B\in \mathcal{B}(\mathcal{P})$, we can write:
\begin{align*}\P_x(R_1\in B)& =\int_0^{+\infty}\P_x(X(u)\in B)e^{-u}\mathrm{d}u
\\ &\geqslant  e^{-N(x)\overline{T}}\int_0^{\overline{T}}\P_x\left[h_x^{N(x)}\left(S^{N(x)}_u,\left(U_1,\ldots,U_{N(x)-1}\right),\kappa_x\right)\in B\right]\mathrm{d}u \end{align*}
where $S^{N(x)}_u=\left(T_1,T_2-T_1,\ldots,T_{N(x)-1}-T_{N(x)-2},u\right)$.\pe

Computing the joint density of the inter-jump times and the $U_k$ and considering a common upper bound $\overline{\theta}$ for the $\theta_{i,j}$ on $\{x'\in \R_+^n \mid \forall y\in C, \Vert x' - y \Vert_1\leqslant N(x)\overline{T}M\}$ yields, for all $x\in C$ and all $B\in \mathcal{B}(\mathcal{P})$:
$$\P_x(R_1\in B)\geqslant e^{-\left(1+\overline{\theta}n(n-1)\right)N(x)\overline{T}}\int_{[0,\overline{T}]^{N(x)}}\int_{[0,1]^{N(x)-1}}
\left[\prod_{i=1}^{N(x)-1}\theta_{\kappa_i(x)}
\left(h_x^i(t,\xi,\kappa_x)\right)\right]\mathrm{1}_{h_x^{N(x)}(t,\xi,\kappa_x)\in B}\mathrm{d}\xi\mathrm{d}t,$$
and therefore
$$\P_x(R_1\in B)\geqslant e^{-\left(1+\overline{\theta}n(n-1)\right)N(x)\overline{T}}\int_{U_2^x}\int_{U_1^x}
\left[\prod_{i=1}^{N(x)-1}\theta_{\kappa_i(x)}
\left(h_x^i(\sigma_x(z_1,z_2),\kappa_x\right)\right]\mathrm{1}_{h_x^{N(x)}(\sigma_x(z_1,z_2),\kappa_x)\in B}\mathrm{d}z_1\mathrm{d}z_2$$
by the change of variables formula. Applying this formula again using Assumption \ref{assump_petite}.\ref{hyp_dilat} and recalling Assumption \ref{assump_petite}.\ref{hyp_mino_dens_sauts} yields
$$\P_x(R_1\in B)\geqslant e^{-\left(1+\overline{\theta}n(n-1)\right)N(x)\overline{T}}\frac{\theta_0}{\lambda_{q(x)}(U_1^x)}\int_{U_2^x}\frac{1}{M}\int_{\psi^x_{z_2}(U_2^x)}\mathrm{1}_{y\in B}\mathrm{d}y\mathrm{d}z_2.$$
Using Assumption \ref{assump_petite}.\ref{constant_scanned}, we may thus write:
$$\P_x(R_1\in B)\geqslant e^{-\left(1+\overline{\theta}n(n-1)\right)N(x)\overline{T}}~\frac{\theta_0}{M}\cdot \lambda_p(B)$$
for any $B\in \mathcal{B}(\mathcal{P})$, which entails that $C$ is petite for $(R_n)_{n\in \N}$, then petite for $X$.
\end{demprop}


\subsection{Ergodicity criterion \ref{thm_ergo}}

\begin{demthm}[\ref{thm_ergo}]
Theorems 7.1 and 7.2 of \cite{MT2} ensure it is sufficient for $F$-ergodicity to hold to prove that there are $\delta>0$ and a compact set $C\subset \R_+^n$ such that
\begin{equation}\label{premier_result_finitude}\sup_{x\in C} ~\E_x\left[\int_0^{\tau_C(\delta)}F\left(X(t)\right)dt\right]<+\infty\end{equation} where $\tau_C(\delta):=\inf\{t\geq \delta\mid X(t)\in C\}$, and that:
\begin{equation}\label{second_result_finitude}\forall x\in \R_+^n,\quad \E_x\left[\int_0^{\tau_C(0)}F\left(X(t)\right)dt\right]<+\infty.\end{equation}
Most of the work we need to provide in order to prove that the assumptions we stated in Section \ref{section_bornitude} entail (\ref{premier_result_finitude}) and (\ref{second_result_finitude}) has already been done in the process of proving Theorem \ref{thm_born_prob}, where we showed that the hitting time $\tau_C$ of any compact subset $C=\{x\in \R_+^n\mid \Vert x \Vert_1\leqslant R\}$ was such that $x\mapsto \E_x\left[\exp\left(\beta \sqrt{\tau_C}\right)\right]<+\infty$ was finite-valued and locally bounded on $\R_+^n$ for some $\beta>0$. \pe 

We now show that there exists $\eta\in \R_+^*$  such that defining $F$ by (\ref{defF}),
$$x\mapsto \E_x\left[\int_0^{\tau_{C}(1)}F\left(X(t)\right)dt\right]$$
is locally bounded on $\R_+^n$. This will end the proof of Theorem \ref{thm_ergo}.\pe 

Using the notations of the proof of Theorem \ref{thm_born_prob}, we set $\eta_0= \frac{\beta}{\sqrt{M}}$ and
$$F_0:\begin{cases}\R_+\to \R_+\\ y\mapsto \mathrm{1}_{y>\frac{1}{\eta_0^2}}\frac{e^{\eta_0 \sqrt{y}}}{\sqrt{y}}\end{cases}$$ 
For any $x\in \R_+^n$, we observe that
\begin{equation}\label{ineg_intlag}\E_x\left[\int_0^{\tau_C(1)}F_0\left(\sum_{i=1}^n X_i(t)\right)\mathrm{d}t\right]\leqslant \E_x\left[\int_0^1 F_0\left(\sum_{i=1}^n X_i(t)\right)\mathrm{d}t+\E_{X(1)}\left[\int_0^{\tau_C}F_0\left(\sum_{i=1}^n X_i(t)\right)\mathrm{d}t\right]\right]\end{equation}
by the Markov property and that for any $X$-stopping time $\tau$ and any $z\in \R_+^n$:
\begin{align*}\E_z\left[\int_0^{\tau}F_0\left(\sum_{i=1}^n X_i(t)\right)\mathrm{d}t\right]
 \leqslant ~&\E_z\left[\int_0^{\tau}F_0\left(\sum_{i=1}^n z_i+tM\right)\mathrm{d}t\right]  
=~ \E_z\left[\frac{2}{\eta_0}\left(e^{\eta_0 \sqrt{\sum_{i=1}^n z_i+\tau M}}-e\right)_+\right] 
& \\ \leqslant ~& \E_z\left[\frac{2}{\eta_0}e^{\eta_0 \sqrt{\sum_{i=1}^n z_i+\tau M}}\right] 
\leqslant~ \frac{2}{\eta_0}e^{\eta_0 \sqrt{\sum_{i=1}^n z_i}}~\E_z\left[e^{\eta_0 \sqrt{\tau M}}\right] 	\end{align*}
since $F_0$ is nondecreasing and $\sqrt{a+b}\leqslant \sqrt{a}+\sqrt{b}$ for all $(a,b)\in \R_+^2$.\pe 

We know from Step 3 of the proof of Theorem \ref{thm_born_prob} that $z\mapsto \E_z\left[e^{\eta_0 \sqrt{\tau_C M}}\right]=\E_z\left[e^{\beta \sqrt{\tau_C}}\right]$ is locally bounded on $\R_+^n$. So are therefore $z\mapsto \E_z\left[\int_0^{\tau_C}F_0\left(\sum_{i=1}^n X_i(t)\right)\mathrm{d}t\right]$ and, consequently, $z\mapsto \E_z\left[\int_0^{\tau_C(1)}F_0\left(\sum_{i=1}^n X_i(t)\right)\mathrm{d}t\right]$ by (\ref{ineg_intlag}) since $\Vert X_1\Vert_1\leqslant \Vert x \Vert_1+M$ holds $\P_x$-almost surely for all $x\in \R_+^n$. We easily find $\alpha_0,\beta_0,\gamma\in \R_+^*$ such that $e^{\gamma \sqrt{y}}\leqslant \alpha_0+\beta_0 F_0(y)$ for all $y\in \R_+$, which entails:
$$\forall x\in \R_+^n,\quad \E_x\left[\int_0^{\tau_C(1)} e^{\gamma \sqrt{\sum_{i=1}^n X_i(t)}}\mathrm{d}t\right]\leqslant \alpha_0~\E_x(\tau_C(1))+\beta_0~\E_x\left[\int_0^{\tau_C(1)}F_0\left(\sum_{i=1}^n X_i(t)\right)\mathrm{d}t\right].$$
This yields the expected result since $x\mapsto \E_x\left[\tau_C(1)\right]$ is locally bounded on $\R_+^n$.
\hfill $\square$\end{demthm}

\subsection{Stability criterion \ref{thm_pcpal_mult} for the multiplicative setting}

We now consider that Assumptions \ref{addass_auto} and \ref{addass_topo} hold and set some further notation.\pe 

First let $\gamma=(\gamma_1,\ldots,\gamma_r)\in \mathcal{E}^r$ be a cycle of $\mathcal{G}$ that visits all sinks. Then, for all $l\in \ie{1}{n}$, let $d_-(l)$ be the $\mathcal{G}$-graph distance from patch $l$ to $V^-$. Given Assumption \ref{addass_topo}, there exists a vector $\kappa=(\kappa_1,\ldots,\kappa_m)\in \mathcal{E}^m$, such that, writing $\kappa_k=(i_k,j_k)$ for all $k\in \ie{1}{m}=V^+\cup V^0$:
\begin{itemize}\setlength\itemsep{0em}
\item For all $k\in V^+\cup V^0$, $d_-(i_k)> d_-(j_k)$ (edges in $\kappa$ are "directed towards the exit of the system", that is, towards sinks). 
\item For all $l\in V^+\cup V^0$, there is a $k(l)\in V^+$ such that $i_{k(l)}=l$ (each source or neutral patch is the origin and one and only one edge in $\kappa$).
\end{itemize}
This choice of $\kappa$ corresponds to an ordered collection of edges along which the population in every patch of the system "reaches the exit by the shortest path". A possible choice of $\kappa$ when $\mathcal{G}$ is the graph in Figure \ref{fig_ex_graphe_admissible} is $\kappa=\left((2,1),(3,1),(4,1),(1,5)\right)$ ; if $\mathcal{G}$ is the first graph of Figure \ref{importance_struc_graphe}, then we may choose $\kappa=\left((1,3),(2,4)\right)$.\pe 

Recalling Proposition \ref{prop_petits} and Theorem \ref{thm_ergo}, it is sufficient to prove the following:

\begin{framed}\vspace*{-.45cm}\begin{prop}\label{prop_thm_mult}
Assume that the model is multiplicative and (\ref{cdt_flux}) holds. Then:
\begin{enumerate}[label=(\roman*)]
\item \label{bornitude_thm_mult}There exists $S\subset \R_+^n$ such that $X$ fits Assumptions \ref{hyp_zone_decr} and \ref{assump_bound}.
\item \label{petitesse_thm_mult}$X$ fits Assumption \ref{assump_petite}.
\item \label{crit_irred_squel_mult}$X$ admits an irreducible skeleton chain.
\end{enumerate}
\end{prop}\vspace*{-.45cm}\end{framed}

\begin{demprop}[\ref{prop_thm_mult}.\ref{bornitude_thm_mult}] Let $q>0$ and $a>0$ be such that $\mu_{i,j}\left([a,1-a]\right)>q$ for all $(i,j)\in \mathcal{I}$. Using (\ref{cdt_flux}), we may choose $R>0$ such that Assumption \ref{hyp_zone_decr} holds for $S=\{x\in \R_+^n\mid \min_i x_i \geqslant R\}$. Let us now set $Z_0=2na^{-(m+r)}\left(M+R\right)$, $M_0=\inf_{i\in V^+}\inf_{[0,Z_0]}\phi_i$, $\underline{\theta}=\min_{(i,j)\in \mathcal{E}}\theta_{i,j}$, $\overline{\theta}=\max_{(i,j)\in \mathcal{E}}\theta_{i,j}$ and
$$S'=\left\{x\in \R_+^n\mid \min_{i\in \ie{1}{n}} x_i \geqslant 2R\right\}.$$
We will show that $S'$ meets Assumption \ref{assump_bound}.\ref{majo_temps_retour_zone_decr} by considering suitable paths of $X$ defined by transfers along the edges of $\kappa$ and $\gamma$.\pe

\indent It is an easy matter to see that with probability at least $\exp\left(-|\mathcal{E}|\left(\frac{Z_0}{M_0}+2\right)\overline{\theta}\right)\underline{\theta}^{m+r}q^{m+r}\frac{1}{m^m r^r}$, the following holds for the path of $X$ stemmed from any $x\in \R_+^n$:
\begin{itemize}\setlength\itemsep{0em}
\item No transfer occurs before time $\frac{Z_0}{M_0}$.
\item $m$ successive transfers occur along edges $\kappa_1,\kappa_2,\ldots$ and $\kappa_m$ between times $\frac{Z_0}{M_0}$ and $\frac{Z_0}{M_0}+1$. Each of these transfers, if originated from patch $i$ and directed towards patch $j$, has an amplitude between $z$ and $\left(1-a\right)z$.
\item $r$ successive transfers happen along edges of $\gamma$ between times $\frac{Z_0}{M_0}+1$ and $\frac{Z_0}{M_0}+2$, the first of which is undertook from the edge of $V^-$ with the largest population. Each of these transfers, if originated from patch $i$ with population $z$ and directed towards patch $j$, has an amplitude between $az$ et $\left(1-a\right)z$.
\item No other transfer than those just described occurs before time $\frac{Z_0}{M_0}+2$.
\end{itemize}
By construction, on such event one has $X_i(1)\geqslant \frac{Z_0}{n}a^{m+r}-2M=2R$ for all $i\in \ie{1}{n}$ so $X(1)\in S'$, which entails that Assumption \ref{assump_bound}.\ref{majo_temps_retour_zone_decr} holds.\pe 

Besides, Assumption \ref{assump_bound}.\ref{mino_temps_depart_zone_decr} is met for any choice of $\varepsilon\in [0,1)$ and $T'>0$ if $R$ is large enough. Indeed, if we define
$$A_{i,j}^{\varepsilon'}(T')=\bigcup_{\substack{k\geqslant 1 \\ T_k<T'}}\left(q_{i,j}(X(T_k^-),U_k)\leqslant (1-\varepsilon')X_i(T_k^-)\right)$$
for any $T'>0$, $(i,j)\in \mathcal{I}$ and $\varepsilon'>0$, $A_{i,j}^{\varepsilon'}(T')$ being the event that all $i\to j$ transfers before time $T'$ are of relative amplitude less than $1-\varepsilon'$, one can write that for any $T'>0$,
$$\P_x\left(A_{i,j}^{\varepsilon'}(T')\right)\underset{\varepsilon'\to 0}{\longrightarrow}1$$
uniformly in $x$ since the $\mu_{i,j}$ assign mass $0$ to $\{1\}$. Now the $\phi_i$ and $\theta_{i,j}$ being bounded implies that for all $T'>0$:
$$\liminf_{\substack{\Vert x \Vert_1 \to +\infty \\ x\in S'}}\P_x(\forall s\in [0,T'],~X(s)\in S)\geqslant \liminf_{\substack{\Vert x \Vert_1 \to +\infty \\ x\in S'}}\P_x\left(\bigcap_{(i,j)\in \mathcal{E}}A_{T'}^{\varepsilon'}(i,j)\right)$$
for all $\varepsilon'\in (0,1)$, hence the result as $\varepsilon'$ tends to $0$. Assumption \ref{assump_bound}.\ref{effic_phases_descente} is then true if $T'$ is large enough (which is possible as soon as $R$ is), which completes our proof of Proposition \ref{prop_thm_mult}.\ref{bornitude_thm_mult}.
\end{demprop}

\begin{demprop}[\ref{prop_thm_mult}.\ref{petitesse_thm_mult}] Using Proposition \ref{transmission_petitesse}, we will consider the behavior of our process starting from a small ball centered on a state that corresponds to positive population levels for sources and neutral patches, and use the change of variables formula to check for Assumption \ref{assump_petite}.\ref{hyp_dilat}.\pe 

\indent It is clear that we may assume that $C=\{x\in \R_+^n\mid \Vert x \Vert_1\leqslant R\}$ for some $R>0$.

\paragraph*{Step 1: reachability of the $\overline{\mathcal{B}}_{\infty}(x^*,\delta)$ for some $x^*$}
~\pe \\
\indent For any $x\in \R_+^m\times \{0\}^{n-m}$ and $\delta>0$, we denote the closed ball of $\R_+^m\times \{0\}^{n-m}$ for the infinite norm by $\overline{\mathcal{B}}_{\infty}(x,\delta)$. We will show that there exists $x^*\in \left(\R_+^*\right)^m\times \{0\}^{n-m}$ such that any $\overline{\mathcal{B}}_{\infty}(x^*,\delta)$ fulfills the requirement of Proposition \ref{transmission_petitesse}.\pe 

Let us first assume that $V^0=\varnothing$. By assumption on $\Phi$, there exists $T>0$ such that $\Phi_i(x,T)=0$ for all $x\in \R_+^n$ with $\Vert x \Vert_1\leqslant R+M$ and all $i\in V^-$. Let $x^*=\Phi\left(0,T\right)$ and first note that $x^*=\Phi\left(z,T\right)$ for all $z\in \R_+^n$ such that $\Vert z \Vert_1\leqslant R+M$ and $z_i=0$ for all $i\in V^0\cup V^-$. Now consider $\delta>0$. $(\Phi_1,\ldots,\Phi_m)$ is uniformly continuous on $C\times [T,T+1]$ since the non-negative $\phi_i$ are $\mathcal{C}^1$, so there are $\delta'\in (0,M]$ and $t_0\in (0,1]$ such that $\Vert \Phi(x,u)-x^* \Vert_{\infty} \leqslant \delta$ for all $x\in \R_+^n$ with $x_1+\ldots+x_m\leqslant \delta'$ and all $u\in [T,T+t_0]$. Set $\varepsilon\in (0,1)$ and $\eta\in \left(0,1\right)$ such that: \begin{equation}\label{varepilon_et_a}\forall (i,j)\in \mathcal{I}, \quad \overline{\mu}_{i,j}([0,\varepsilon))< \eta\end{equation}
and denote by $N$ the smallest positive integer such that $\left(1-\varepsilon\right)^NR<\frac{\delta'}{2}$.\\

It is easy to see that setting $\kappa^0=(\kappa,\ldots,\kappa)\in \mathcal{I}^{mN}$, we have:
$$\forall x\in C,\quad  \forall t\in \left[0,\frac{\delta'}{2mNM}\right]^{mN}\times \left[T,T+t_0\right],\quad  \forall \xi\in \left(\eta,1\right]^{mN},\quad h^{mN+1}_x(t,\xi,\kappa^0)\in \overline{\mathcal{B}}_{\infty}\left(x^*,\delta\right).$$ 
Considering the joint density of the $T_k$ and the $U_k$ as in the proof of Proposition \ref{prop_thm_mult}.\ref{bornitude_thm_mult} finally yields, for all $x\in C$:
\begin{equation}\label{ineg_petitesse_multipl}\P_x\left(\forall t\in \left[T,T+\frac{\delta'}{2M}+t_0\right], X(t)\in \overline{\mathcal{B}}_{\infty}\left(x^*,\delta\right)\right)\geqslant \exp\left(-|\mathcal{E}|\left(\frac{\delta'}{2M}+t_0\right)\overline{\theta}\right)\left(\frac{\underline{\theta}\delta'}{2mNM}\right)^{mN}(1-\eta)^{mN}\end{equation}
which, in turn, implies the expected property.\pe\pe

The proof relies on the same reasoning if $V^0\neq \varnothing$, using paths along which a sequence of transfers results in all neutral patches having a positive population, then the flow brings the population of sinks to zero.

\paragraph*{Step 2: petiteness of a $\overline{\mathcal{B}}_{\infty}(x^*,\delta)$ for some $\delta$}
~\pe \\
\indent According to Proposition \ref{transmission_petitesse} and Step 1 hereabove, Proposition \ref{prop_thm_mult}.\ref{petitesse_thm_mult} will be proved if we show that there exists some $\delta>0$ such that $\overline{\mathcal{B}}_{\infty}(x^*,\delta)$ is petite for the resolvent of $X$, which can be done by way of checking Assumption \ref{assump_petite}.\pe 

Let us thus set some $\delta\in (0,\min_{i\in V^+\cup V^0} x^*_i)$ that we may have to take smaller later on, and let $\underline{x}=\min_{i\in V^+\cup V^0}x_i^*-\delta$.\pe

For simplicity reasons, we will assume that all $\phi_i$ associated with sources are $\mathcal{C}^1$ functions (rather than merely piecewise $\mathcal{C}^1$) on $\R_+$ and that $\overline{\mu}_{i,j}$ measures admit continuous density functions $f_{i,j}$ on $[0,1]$ (instead of just admitting an absolutely continuous component on a subinterval of $[0,1]$). The general case only requires reducing the domain over which it is possible to consider our paths of interest, which induces an unnecessary notational inflation.\pe 

Even if it means considering a larger $T$ in Step 1, we may assume that $\Phi(x,T)\in \R_+^{m}\times \{0\}^{n-m}$ for all $x\in \overline{\mathcal{B}}_{\infty}(x^*,\delta)$. As in Step 1, we define $\varepsilon$ and $\eta$ such that (\ref{varepilon_et_a}) holds.\pe 

Let us begin by proving the following inequality: 
\begin{equation}\label{mino_der_flux}\forall i\in V^+,\forall y>0, \forall u>0, \quad \frac{\min_{[0,y+Mu]}\phi_i}{M}\leqslant  \frac{\partial \Phi_i}{\partial y}(y,u) \leqslant \frac{M}{\min_{[0,y+Mu]}\phi_i}.\end{equation}
If $i\in V^+$, $0<y<y+h$ and $u>0$, the mean value inequality entails that $\Phi_i\left(y,\frac{h}{M}\right)\leqslant y+h$, from which we deduce
$$\frac{\min_{\left[\Phi_i(y,u),\phi_i(y,u)+h\right]}\phi_i}{M}h \leqslant \phi_i\left(y,u+\frac{h}{M}\right)-\Phi_i\left(y,u\right)\leqslant \Phi_i(y+h,u)-\Phi_i\left(y,u\right)$$
using the mean value theorem, hence the left hand side of (\ref{mino_der_flux}) by letting $h$ go to $0$. Proving the right-side inequality of (\ref{mino_der_flux}) relies on the very same argument and is left to the reader. Moreover, it is clear that \begin{equation}\label{mino_der_flux_nul}\forall i\in V^0,\forall y>0, \forall u>0, \quad \frac{\partial \Phi_i}{\partial y}(y,u)=1.\end{equation}

We now show that for any $x\in \overline{\mathcal{B}}_{\infty}(x^*,\delta)$ and $t\in (0,1)^m\times (T,T+1)$, then
$$\psi^x_{t}:\begin{cases} (0,a)^m \to \R^m \times \{0\}^{n-m} \\ \xi\mapsto h_x^{m+1}\big(t,\xi,\kappa\big)\end{cases}$$
is a $\mathcal{C}^1$-diffeomorphism of $(0,\eta)^m$ onto its image. For fixed $x$ and $t$, indeed, a simple calculation shows that 
\begin{equation}\label{der_part_1}\frac{\partial \psi^x_{t}}{\partial \xi_m}(t,\xi,\kappa)=\partial_1 \Phi\left(g_{\kappa_m}(h^m_x,\xi_m),t_{m+1}\right)\cdot {\partial}_2 g_{\kappa_m}(h^m_x,\xi_m)\end{equation}
and that for any $i\in \ie{1}{m-1}$:
\begin{equation}\label{der_part_2}\frac{\partial \psi_x^{t}}{\partial \xi_i}(t,\xi,\kappa)= \left[\prod_{k=m}^{i+1} \partial_1 \Phi\left(g_{\kappa_k}(h^k_x,\xi_k),t_{k+1}\right)\cdot {\partial}_1 g_{\kappa_k}\left(h^k_x,\xi_k\right)\right]\cdot \partial_1 \Phi\left(g_{\kappa_i}(h^i_x,\xi_i),t_{i+1}\right) \cdot {\partial}_2 g_{\kappa_i}(h^i_x,\xi_i).\end{equation}
where we wrote $h_x^k$ as a short for $h_x^k(t,\xi,\kappa)$, defined
$${\partial}_1 g_{i,j}(y,\zeta)=I_m+F_{\mu_{i,j}}^{-1}(\zeta)\cdot(\mathrm{1}_{j\leqslant m}E_{j,i}-E_{i,i})$$
with $\left(E_{i,j}\right)_{i,j\in \ie{1}{m}}$ the canonical basis of $\mathcal{M}_m(\R)$, as well as
$${\partial}_1\Phi(y,u)=\begin{pmatrix}\frac{\partial\Phi_1(y,u)}{\partial y_1} & 0 & \cdots & 0 \\ 0 & \ddots & \ddots & \vdots \\ \vdots &  \ddots & \ddots & 0 \\ 0 & \cdots & 0 & \frac{\partial\Phi_m(y,u)}{\partial y_m}\end{pmatrix}$$
and
$${\partial}_2 g_{i,j}(y,\xi)=\frac{y_i}{f_{i,j}\left(F^{-1}_{\mu_{i,j}}\left(\xi\right)\right)}\left(\mathrm{1}_{j\leqslant m}\tilde{e}_j-\tilde{e}_i\right)$$
where $\left(\tilde{e}_1,\ldots,\tilde{e}_m\right)$ is the canonical basis of $\R^m$.\pe 

Recalling the fact that every $i\in V^+$ is the origin of one only edge $\kappa_i$,  (\ref{mino_der_flux}) and (\ref{mino_der_flux_nul}) then imply that for all $t\in (0,1)^m\times (T,T+1)$ and all $\xi\in (0,\eta)$, the Jacobian matrix of $\psi^x_t$ at $(t,\xi,\kappa)$ is a continuous function of $\xi$ and is an invertible matrix with determinant $J(t,\xi,\kappa)$ such that
\begin{equation}\label{def_alpha_mult}0 < J(t,\xi,\kappa)\leqslant \left(\frac{M}{\min_{i\in V^+}\min_{[0,x_i+(m+T+1)M]}\phi_i}\right)^{\frac{m(m+1)}{2}}\frac{\left(\Vert x \Vert_1 + \delta + (m+T+1)M\right)^m}{\prod_{i=1}^m \min_{[0,1]}f_{\kappa_i}} <+\infty.\end{equation}
Now $\psi^x_t$ is injective over $(0,\eta)^m$ for fixed $t$ and $x$; this stems from the fact that the cumulative distribution functions of the $\mu_{\kappa_k}$ are strictly increasing, so knowing $\langle  \psi^x_t(t,\xi,\kappa), e_{i_1}\rangle$ makes it possible to determine $\xi_1$, knowing $\xi_1$ and $\langle  \psi^x_t(t,\xi,\kappa), e_{i_2}\rangle$ yields $\xi_2$ and so on. This proves our claim that $\psi^x_t$ is a $\mathcal{C}^1$-diffeomorphism.\pe

We therefore only need to check that Assumption \ref{assump_petite}.\ref{constant_scanned} holds to conclude. Take $x\in \overline{\mathcal{B}}_{\infty}(x^*,\delta)$ and for all $t\in [0,1]^m$, set $I_t=\psi^x_t(t,(0,\eta)^m,\kappa)$. Considering the possible value of the population of patches $1$ to $m$ after a series of jumps along the edges of $\kappa$ such that the $k$-th transfer has amplitude between $0$ and $k\frac{\varepsilon}{m} \underline{x}$, it is fairly easy to check that for a given $t=(t_1,\ldots,t_m)$ in $\left(0,\left(\frac{\varepsilon}{Mm^2}\underline{x}\right)\wedge 1\right)^m\times (T,T+1)$:
$$\prod_{i=1}^m \left[\Phi_i\left(x_i-\frac{\varepsilon}{m}+(t_1+\ldots+t_i)M,t_{i+1}+\ldots+t_{m+1}\right),\Phi_i\left(x_i,t_1+\ldots+t_{m+1}\right)\right)\subset I_t$$
and then, by (\ref{mino_der_flux}), (\ref{mino_der_flux_nul}) and the mean value theorem for the $\Phi_i(\cdot,t)$:
$$\prod_{i=1}^m \left[\Phi_i(x_i,t_{i+1}+\ldots+t_{m+1})-\varepsilon_0,\Phi_i(x_i,t_1+\ldots+t_{m+1})\right)\subset I_t$$
for some positive $\varepsilon_0$ independent from the choice of $x$ within $\overline{\mathcal{B}}_{\infty}\left(x^*,\delta\right)$. For small enough values of $\delta$, there are $T^0\in \left(0,\left(\frac{\varepsilon}{Mm^2}\underline{x}\right)\wedge 1\right]$ and an orthotope $\mathcal{P}\subset \R_+^m\times \{0\}^{n-m}$ with non-zero Lebesgue measure, both independent from the choice of $x$ within $\overline{\mathcal{B}}_{\infty}\left(x^*,\delta\right)$, such that $\mathcal{P}\subset I_t$ for all $t\in \left(0,T^0\right)^m\times (T,T+T^0)$. This entails that Assumption \ref{assump_petite}.\ref{constant_scanned} holds. Proposition \ref{prop_petits} then implies that $\overline{\mathcal{B}}_{\infty}(x^*,\delta)$ is petite for the resolvent of $X$, which ends the proof of Proposition \ref{prop_thm_mult}.\ref{petitesse_thm_mult}, and therefore of Theorem \ref{thm_pcpal_mult}.\hfill $\square$\pe\enlargethispage{4mm}

Note that if the $t\mapsto h_x^m(t,\xi,\kappa)$ are diffeomorphisms for suitable values of $\xi$ (which typically is the case in a multiplicative model with constant growth!), we may want to monitor the effect of small variations of $t$ (rather than of $\xi$) on $h_x^m(t,\xi,\kappa)$. This would allow to relax the absolute continuity assumptions on the $\mu_{i,j}$, since it would then be sufficient to assume that $\overline{\mu}_{i,j}\left(\{0\}\right)<1$ and $\overline{\mu}_{i,j}\left(\{1\}\right)=0$ for petiteness to hold.\end{demprop}

\begin{demprop}[\ref{prop_thm_mult}.\ref{crit_irred_squel_mult}]
We keep all notations from the proof of Proposition \ref{prop_thm_mult}.\ref{petitesse_thm_mult} above. First set $\Delta_0 = \frac{\delta'}{2M}+t_0$. If $\pi$ denotes the invariant probability of $X$, we may assume that $\pi(C)>0$, which entails by Birkhoff's ergodic theorem that $\tau_C$ is $\P_x$-a.s. finite for all $x\in \R_+^n$. Inequality (\ref{ineg_petitesse_multipl}) then yields:
$$\inf_{x\in \R_+^n}\P_x\left(\forall t\in \left[\tau_C+T,\tau_C+T+\Delta_0\right], X(t)\in \overline{\mathcal{B}}_{\infty}\left(x^*,\delta\right)\right)>0.$$
Thus, according to the strong Markov property, it is sufficient to show that there exists $\Delta\in (0,\Delta_0]$ such that: 
$$\forall x\in \mathcal{B}(x^*,\delta), \forall B\in \mathcal{B}(\mathcal{P}),\quad \lambda_m(B)>0\Rightarrow \P_x\left(\exists k\geqslant 1, X(k\Delta)\in B\right)>0.$$ 
Set $\Delta=\Delta_0\wedge T$ and $k\in \Z_+$ such that $mT^0\leqslant k\Delta<mT^0+T$, and let $x\in \overline{\mathcal{B}}_{\infty}(x^*,\delta)$. Then for any Borel subset $B$ of $\mathcal{P}$, the change of variables formula yields:
\begin{align*}\P_x\left(X(k\Delta)\in B\right) ~\geqslant~ & \P_x\left[\left(T_m<k\Delta<T_{m+1}\right)\cap\left(\Phi(X(T_m),k\Delta-{T_m})\in B\right)\right] \\
~\geqslant~ & e^{-\overline{\theta}n(n-1)k\Delta}\int_{(0,T^0)^{m}}\int_{(0,\eta)^m}\left[\prod_{i=1}^{m}\theta_{\kappa_i(x)}
\left(h_x^i\right)\right]\mathrm{1}_{\Phi\left(h_x^{m},k\Delta-\sum_{j=1}^{m}t_j\right)\in B}\mathrm{d}\xi\mathrm{d}t & \\
 \geqslant~ & e^{-\overline{\theta}n(n-1)k\Delta} \underline{\theta}^m ~\int_{(0,T^0)^m}\int_{(0,\eta)^m} \mathrm{1}_{\Phi\left(h_x^{m},k\Delta-\sum_{j=1}^{m}t_j\right)\in B}\mathrm{d}\xi\mathrm{d}t & \\
 \geqslant~ & e^{-\overline{\theta}n(n-1)k\Delta} \underline{\theta}^m ~\int_{(0,T^0)^m}\int_{(0,\eta)^m} \mathrm{1}_{h_x^{m+1}\left((t,k\Delta-\sum_{j=1}^{m}t_j),\xi,\kappa\right)\in B}\mathrm{d}\xi\mathrm{d}t & \\ 
 \geqslant~ & \frac{1}{\alpha}e^{-\overline{\theta}n(n-1)k\Delta} \underline{\theta}^m ~\int_{(0,T^0)^m}\mathrm{1}_{k\Delta-\sum_{j=1}^{m}t_j\in [T,T+T^0]}\left(\int_{\mathcal{P}} \mathrm{1}_{y\in B}\mathrm{d}y\right)\mathrm{d}t & \\
 \geqslant~ & \frac{1}{\alpha}e^{-\overline{\theta}n(n-1)k\Delta} \underline{\theta}^m  \lambda_m\left(S\right)\lambda_m(B)
\end{align*}          
where $h_x^i$ stands for $h_x^i(t,\xi,\kappa)$, $\alpha$ is the upper bound of $J(t,\xi,\kappa)$ defined by (\ref{def_alpha_mult}) and
$$S=\left\{t\in (0,T^0)^m, \sum_{j=1}^m t_j \in \left[k\Delta-T-T^0,k\Delta-T\right]\right\}.$$
Now $\lambda_m\left(S\right)>0$, which ends the proof.\hfill $\square$\end{demprop}

\subsection{Stability criterion \ref{thm_pcpal_uni} for the unitary setting}

We keep considering $\gamma$ and $\kappa$ as defined in the proof of Theorem \ref{thm_pcpal_mult}. We will show the following, which is sufficient to conclude thanks to Proposition \ref{prop_petits} and Theorem \ref{thm_ergo}.

\begin{framed}\vspace*{-.45cm}\begin{prop}\label{prop_thm_uni}
Assume that the model is unitary and that either (\ref{flux_strong}) holds, or $\mathcal{G}$ is strongly connected and (\ref{flux_weak}) holds. Then:
\begin{enumerate}[label=(\roman*)]\setlength\itemsep{0em}
\item \label{bornitude_thm_uni}There exists $S\subset \R_+^n$ such that $X$ fits Assumptions \ref{hyp_zone_decr} and \ref{assump_bound}.
\item \label{petitesse_thm_uni}$X$ fits Assumption \ref{assump_petite}.
\item \label{crit_irred_squel_uni}$X$ admits an irreducible skeleton chain.
\end{enumerate}
\end{prop}\vspace*{-.45cm}\end{framed}

\begin{demprop}[\ref{prop_thm_uni}.\ref{bornitude_thm_uni}]
Let us first assume that (\ref{flux_strong}) holds. Let $R\geqslant 2$ be an integer such that Assumption \ref{hyp_zone_decr} holds with
$$S=\left\{x\in \R_+^n\mid \min_{i\in V^-}x_i\geqslant R\right\}$$ and set
$$S'=\left\{x\in \R_+^n\mid \min_{i\in V^-}x_i\geqslant 2R\right\}.$$
In order to verify that Assumption \ref{assump_bound}.\ref{majo_temps_retour_zone_decr} holds, we consider $x\in \R_+^n$ as well as $T\in (0,1]$ and argue as in the proof of Proposition \ref{prop_thm_mult}.\ref{bornitude_thm_mult}. On the event we consider, transfers along edges $\kappa_1, \ldots, \kappa_m$ occur before time $T$ and result in the total population of sinks being at least $2(n-m)R+TM$, and transfers along $\gamma$ occur between time $T$ and time $2T$ and result in the population of each sink at time $2T$ reaching at least $2R$. For corresponding paths, $X(2T)\in S'$, and it is easy to show using the assumptions on $\Theta$ that the probability under $\P_x$ of observing such paths for a fixed $T$ goes to $1$ as $\Vert x \Vert_1$ goes to $+\infty$. As a result, $T$ may be taken arbitrarily small in Assumption \ref{assump_bound}.\ref{majo_temps_retour_zone_decr} for any fixed value of $\delta$. Assumption \ref{assump_bound}.\ref{effic_phases_descente} will thus be automatically fulfilled for some value of $T$ provided that Assumption \ref{assump_bound}.\ref{mino_temps_depart_zone_decr} holds. Showing the latter relies on the observation that as a sink gives up its charge, the temporal intensity of subsequent transfers from this patch is upper bounded. Simple calculations then show that the probability for a sink with original population above $2R'$ to have a population lower than $R'$ before time $\frac{R'}{M}$ is lower than some constant in $[0,1)$, which ends our proof.\pe

The argument is similar in the connected case with (\ref{flux_weak}) except that it is now possible to define $S$ and $S'$ as in the proof of Proposition \ref{prop_thm_mult}.\ref{bornitude_thm_mult}. The connectivity assumption  then entails that the process returns to $S'$ arbitrarily fast with a given probability when the total population in the system is high enough, and we may conclude just as before.\enlargethispage{4mm}
\end{demprop}

\begin{demprop}[\ref{prop_thm_uni}.\ref{petitesse_thm_uni}]
Without loss of generality, we may assume that $C$ writes $\{x\in \R_+^n\mid \Vert x \Vert_1\leqslant R\}$ for some $R\in \Z_+^*$. Let $T>0$ be such that $\Phi_i(x,T)=0$ for all $i\in V^-$ and all $x\in \R_+^n$ with $\Vert x \Vert_1\leqslant R+1$ and define
$$\kappa'=\left(\kappa_1,\ldots,\kappa_1,\kappa_2,\ldots,\kappa_2,\ldots,\kappa_m,\ldots, \kappa_m\right)\in \mathcal{I}^{m(R+1)}$$
as the vector of edges made of $R+1$ successive repetitions of each of $\kappa$'s edges. 
\pe \\
For all $t\in \left[0,\frac{1}{mM(R+1)}\right)^{m(R+1)}\times (T,T+1)$ and all $\xi\in (0,1]^{m(R+1)}$, we may write:
$$\forall i\in V^+, \quad \left(h_x\left(t,\xi,\kappa'\right)\right)_i=\Phi_i\left(0,t_{(i+1)(R+1)}+t_{(i+2)(R+1)}+\ldots+t_{m(R+1)}+t_{m(R+1)+1}\right)$$
and
$$\forall i\in V^0\cup V^-, \quad \left(h_x\left(t,\xi,\kappa'\right)\right)_i=0.$$
For $t'=(t_1,\ldots,t_{R},t_{R+2},\ldots,t_{2(R+1)-1},t_{2(R+1)+1},\ldots,t_{m(R+1)-1},t_{m(R+1)+1})\in \left[0,\frac{1}{mM(R+1)}\right)^{mR}\times {(T,T+1)}$ and $\xi\in (0,1]^{mR}$, we easily see that
$$\psi^x_{t'}:\begin{cases} \left(0,\frac{1}{mM(R+1)}\right)^{m} \to \R^d\times \R^{n-d} \\ \left(t_{R+1},t_{2(R+1)}\ldots,t_{m(R+1)}\right)\mapsto h_x^{m(R+1)+1}(t,\xi,\kappa)\end{cases}$$
is a $\mathcal{C}^1$-diffeomorphism of $\left(0,\frac{1}{mM(R+1)}\right)^{m}$ onto its image. It is independent from the choice of $x$ in $C$ and has its Jacobian determinant upper bounded by $\alpha=\prod_{i=1}^d \max_{[0,R+1+(T+1)M]}\phi_i$. It is then possible to define, as in the proof of Proposition \ref{prop_thm_mult}.\ref{petitesse_thm_mult}, $T^0\in \left(0,\frac{1}{mM(R+1)}\right)$ and an orthotope $\mathcal{P}\subset \R_+^d\times \{0\}^{n-d}$ with non-zero Lebesgue measure, independently from the choice of $x\in C$, such that $\mathcal{P}$ be included in the image of $\left(0,T^0\right)^{m}$ by $\psi^x_{t'}$ for all $t'\in \left(0,T^0\right)^{mR}\times (T,T+T^0)$. Assumption \ref{assump_petite}.\ref{hyp_mino_dens_sauts} is also met since $\Theta$ is lower bounded by a positive real number, and Proposition \ref{prop_petits} therefore implies that $C$ is petite.
\end{demprop}

\begin{demprop}[\ref{prop_thm_uni}.\ref{crit_irred_squel_uni}]
We use the notations of the proof of Proposition \ref{prop_thm_uni}.\ref{petitesse_thm_uni} above. It is sufficient to prove the irreducibility property for initial conditions within $C$ because of the same argument as in the proof of Proposition \ref{prop_thm_mult}.\ref{crit_irred_squel_mult}. Even if it means considering smaller $T^0$ and $\mathcal{P}$, one may assume that:
$$\forall t'\in \left(0,\frac{T^0}{2m(R+1)}\right)^{mR}\times \left(0,T^0\right), \quad \mathcal{P}\subset \psi_{t'}^x\left(0,\frac{T^0}{2m(R+1)}\right)^m.$$
Let us consider $\Delta\in \left(0,\frac{T^0}{2}\right)$ and $k\geqslant 1$ such that $k\Delta\in \left(T+\frac{T^0}{2},T+T^0\right)$. Reproducing the calculations of the proof of Proposition \ref{prop_thm_mult}.\ref{crit_irred_squel_mult} yields for any $x\in C$ and any $B\in \mathcal{B}(\mathcal{P})$:
\begin{align*}\P_x\left(X_{k\Delta}\in B\right) ~\geqslant~ & e^{-k\Delta\overline{\theta}}\Theta(0)^{m(R+1)}\int_{\left(0,T^0\right)^{m(R+1)}}\mathrm{1}_{h_x^{m(R+1)}((t,k\Delta-T-\sum_{j=1}^{m(R+1)}t_j),\xi,\kappa')\in B}\mathrm{d}t\\
\geqslant~ & \frac{1}{\alpha}e^{-k\Delta\overline{\theta}}\Theta(0)^{m(R+1)}\int_{\left(0,T^0\right)^{m(R+1)}}\int_{\left(0,\frac{T^0}{2m(R+1)}\right)^{mR}}\lambda_d(\mathcal{P}\cap B)\mathrm{d}t'\\
\geqslant~ & \frac{1}{\alpha}e^{-k\Delta\overline{\theta}}\Theta(0)^{m(R+1)}\left(\frac{T^0}{2m(R+1)}\right)^{mR}\lambda_d(B)
\end{align*}
where $\overline{\theta}$ is a common upper bound for the $\theta_{i,j}$ on $\{y\in \R_+^n\mid \Vert y \Vert_1\leqslant R+k\Delta M\}$. Therefore we have $\P_x\left(X(k\Delta)\in B\right)>0$ whenever $\lambda_d(B)>0$, which ends our proof.\hfill $\square$
\end{demprop}

\section*{Appendix A: On the invariant probability of $X$}

\indent \indent It is clear by the proof of Proposition \ref{prop_petits} that if Assumption \ref{assump_petite} holds and if $X$ admits an invariant probability $\pi$, then $\pi$ dominates the Lebesgue measure on some open subset of an affine subspace of $\R^n$. Its support is yet not necessarily included in this subspace, as can be seen in Figure \ref{fig_illus_loi_lim} below. Moreover, determining conditions for the absolute continuity of $\pi$ restricted to given areas of $\R_+^n$ is a non-trivial matter that still has to be discussed. Interested readers are referred to \cite{Loc} for a discussion about the absolute continuity of $\pi$ in a \emph{house of cards} model that has common specifications with our setting.

\begin{figure}
\begin{center}
\includegraphics[scale=.7]{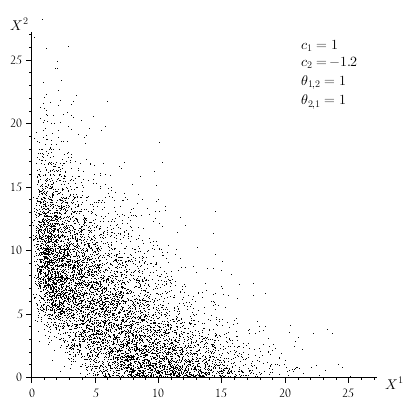}
\includegraphics[scale=.7]{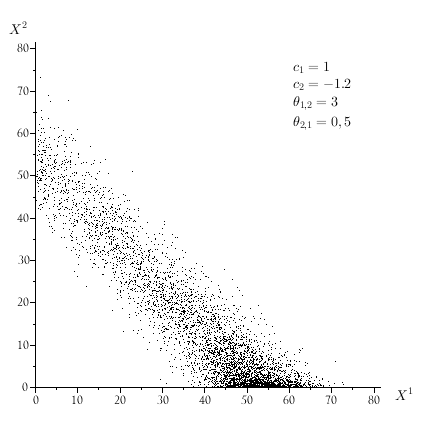}
\captionsetup{justification=justified}
\caption{\label{fig_illus_loi_lim}Plotting of $10000$ simulated instances of $X_{100}$ in the two-patch uniform model with constant growth for specified parameters, with initial value $(5,5)$. In accordance with equation (\ref{busy_time}) of example \ref{exemple_busytime}, about $17.67\%$ of points in both figures lie on the $x$ axis.}
\end{center}
\end{figure}
\indent \indent If $X$ is positive Harris recurrent with invariant probability $\pi$ and $\int \Vert x\Vert_1\mathrm{d}\pi(x)<+\infty$, we see by (\ref{expression_gen}) that for any $\mathcal{C}^1$ function $f:\R_+^n\to \R$ with bounded differential the classical relation hereafter holds:
\begin{equation}\label{gen_annul}\int \left(\sum_{i=1}^n \frac{\partial f}{\partial x_i}(x) \phi^i(x)+\sum_{(i,j)\in \mathcal{I}} \theta_{i,j}(x)\int \left(f(x+\xi(e_j-e_i))-f(x)\right)\mu_{i,j}(x,\mathrm{d}\xi_i)\right)\mathrm{d}\pi(x)=0.\end{equation}
It is worth noting that the second term of the integrand in the LHS above is zero whenever $f(x)$ only depends on $x_1+\ldots+x_n$. Equation (\ref{gen_annul}) then takes a simple form that may provide valuable information on $\pi$. \\
\indent In particular, applying (\ref{gen_annul}) to the projection functions $\langle \cdot,e_i\rangle$ yields:
\begin{equation}\label{eq_coord_sous_pi}\forall t\geqslant 0, \forall i\in \ie{1}{n},\quad \E_{\pi}\left(\phi_i(X(t))\right)+\sum_{j\neq i}\E_{\pi}\left(d_{j,i}(X(t))-d_{i,j}(X(t))\right)=0\end{equation}
where $d_{i,j}:\R_+^n\to \R$ is the \emph{debit function} defined by
\begin{equation}\label{def_debit}d_{i,j}(x)=\theta_{i,j}\left(x\right)\int \xi~\mu_{i,j}(x,\mathrm{d}\xi).\end{equation}
\begin{exemple}\label{exemple_avec_mod_ergoexp}Let us assume that the following conditions hold:
$$\begin{cases}
\phi_i(x)=a_i-x_i\\
\theta_{i,j}(x)=1 \\
\int \xi\mu_{i,j}(x,\mathrm{d}\xi)=m_i x_i
\end{cases}$$ 
where the $a_i$ and $m_i$ are non-negative real numbers, and that $X$ is positive Harris recurrent and integrable under its invariant probability. Then (\ref{eq_coord_sous_pi}) becomes:
$$\forall t\geqslant 0,\forall i\in \ie{1}{n},\quad \left(1+m_i(n-1)\right)\E_{\pi}\left(X_i(t)\right)-\sum_{j\neq i}m_j\E_{\pi}\left(X_j(t)\right)=a_i$$
and can be written as a dominant diagonal linear system, which makes it possible to determine the $\E_{\pi}\left(X_i(t)\right)$ explicitly.
\end{exemple}

\begin{exemple}\label{exemple_busytime}
Let us consider the one-exit constant growth framework, that is, a constant growth model with $c_i\geqslant 0$ for all $i\in \ie{1}{n-1}$ and $c_n<0$, and assume that $X$ is positive Harris recurrent with invariant probability $\pi$ and $\int \Vert x\Vert_1\mathrm{d}\pi(x)<+\infty$. Proving that
\begin{equation}\label{busy_time}\P_{\pi}\left(X_{n}(t)>0\right)=\frac{\sum_{j=1}^{n-1}c_j}{\left|c_n\right|}\end{equation}
and that the event $(X_{n}(t)>0)$ is independent from the variable $\sum_{i=1}^n X_i(t)$ for all $t\geqslant 0$ under $\P_{\pi}$ is left to the reader as an (easy) exercise.
\end{exemple}

\section*{Appendix B: Convergence of a scaled transient multiplicative uniform model towards a Beta distribution}

\begin{demprop}[\ref{gamma_transient}]
We remind that the Beta distribution with parameters $a,b\in \R_+^*$, denoted $\beta(a,b)$, is the probability measure on $\R$ with density
$$f_{a,b}:t\mapsto \mathrm{1}_{]0,1[}(t)\frac{t^{a-1}(1-t)^{b-1}}{B(a,b)}$$ with respect to the Lebesgue measure, where $B(a,b):=\int_0^1 t^{a-1}(1-t)^{b-1}\mathrm{d}t$.\pe 

Set $x\in \mathcal{S}$ and let $S$ be a weak solution of 

\begin{equation*}\mathrm{d}S_1(t)=\int_0^1 \xi S_2(t) N'_{2,1}(\mathrm{d}t,\mathrm{d}\xi)-\int_0^1 \xi S_1(t) N'_{1,2}(\mathrm{d}t,\mathrm{d}\xi)\end{equation*}
with initial value $\frac{x}{\Vert x\Vert_1}$ and such that $S_1+S_2=1$. Here $\overline{\mu}_{1,2}=\overline{\mu}_{2,1}$ is the uniform distribution on $[0,1]$, so $N_{1,2}$ and $N_{2,1}$ have intensity $\mathrm{d}t\mathrm{d}\xi$. Our aim is to show that $S_1(t)$ weakly converges to a $\beta\left(\frac{\theta_{2,1}}{\lambda},\frac{\theta_{1,2}}{\lambda}\right)$ distribution as $t$ tends to infinity. It is easy to see that $S_1$ fulfills the requirements of Theorem 5.2 (c) of \cite{DMT95}, all compact subsets of $\R$ being petite for this process, so we only need to characterize its invariant (and limiting) probability distribution $\pi$. It it clear that the support of $\pi$ lies within $[0,1]$ and that if $f:[0,1]\to \R$ is continuous, then for all $t>0$ the following holds:
$$\int_0^1 \frak{A}'f(s)\mathrm{d}\pi(s)=0$$
where $\frak{A}'$ is the infinitesimal generator of $S_1$, that is:
$$\int_0^1\left[\theta_{1,2}\int_0^1 \left(f(\xi s)-f(s)\right) \mathrm{d}\xi + \theta_{2,1}\int_0^1 \left(f\big(\xi s + (1-\xi)\big)-f(s)\right)\mathrm{d}\xi\right]\mathrm{d}\pi(s) = 0.$$
Applying this equation to power functions $f:x\mapsto x^k$ we get the following recursive formula for the moments of $\pi$:
$$\forall k\geq 1, \quad \int_0^1 s^k \mathrm{d}\pi(s) = \frac{\theta_{2,1}}{k\lambda}\sum_{i=0}^{k-1}\int_0^1 s^i \mathrm{d}\pi(s)$$
which yields, for all $k\geq 1$:
\begin{equation}\label{moments_beta}\int_0^1 s^k \mathrm{d}{\pi}(s)=\frac{1}{k!}\left[ \prod_{i=0}^{k-1}\left(i+\frac{\theta_{2,1}}{\lambda}\right)\right].\end{equation}
Now $\pi$ can be characterized by its moments among Borel probability measures on $[0,1]$ thanks to the Stone-Weierstrass theorem. It is therefore sufficient to compute the moments of the $\beta\left(\frac{\theta_{2,1}}{\lambda},\frac{\theta_{1,2}}{\lambda}\right)$ distribution by using the classical formula  $$B(a,b)=\frac{\Gamma(a)\Gamma(b)}{\Gamma(a+b)}$$ holding for all positive $a$ and $b$, and see that these moments coincide with these given by Equation \ref{moments_beta}.
\hfill $\square$\end{demprop}
%
\bibliography{Biblio}
\bibliographystyle{plain}

\end{document}

%% file: macropulko.tex
\usepackage[T1]{fontenc}
\usepackage{lmodern}
\usepackage{graphicx}
\usepackage{stmaryrd} 
\usepackage{epstopdf}
\usepackage{amsfonts}
\usepackage{dsfont}
\usepackage{amsmath}
\usepackage{array}
\usepackage{amssymb}
\usepackage{enumerate}
\usepackage{euscript}
\usepackage{ulem} 
\usepackage{manfnt}
\usepackage[english]{babel}

\newcommand{\ppe}{\vspace*{0.1cm}}
\newcommand{\pe}{\vspace*{0.2cm}}

\newcommand{\bea}{\begin{array}}
\newcommand{\ena}{\end{array}}

\def\[ent{[\hskip -1.5pt [}
\def\]ent{]\hskip -1.5pt ]}

\newcommand{\ie}[2]{\,\llbracket #1,\,#2 \rrbracket} 

\newcommand \R {\mathbb{R}}

\newcommand \N {\mathbb{N}}
\newcommand \Z {\mathbb{Z}}

\renewcommand{\P}{\mathbb{P}}
\newcommand \E {\mathbb{E}}

%% file: Paperling_revu.bbl
\begin{thebibliography}{10}

\bibitem{ADR67}
J.~Az{\'e}ma, M.~Duflo, and D.~Revuz.
\newblock Mesure invariante sur les classes récurrentes des processus de
  markov.
\newblock {\em Zeitschrift für Wahrscheinlichkeitstheorie und Verwandte
  Gebiete}, 1967.

\bibitem{TCP}
J.-B. Bardet, A.~Christen, A.~Guillin, F.~Malrieu, and P.-A. Zitt.
\newblock Total variation estimates for the tcp process.
\newblock {\em Electronic Journal of Probability}, 18, No 10:1--21, 2013.

\bibitem{Bil99}
Patrick Billingsley.
\newblock {\em Convergence of Probability Measures - Second Edition}.
\newblock Wiley Series in Probability and Statistics, 1999.

\bibitem{BPR}
E.~Brooks-Pollock and M.~J. Roberts, G. O. \&~Keeling.
\newblock A dynamic model of bovine tuberculosis spread and control in great
  britain.
\newblock {\em Nature}, 511:228--231, 2014.

\bibitem{Cas}
C.~G. Cassandras and J.~Lygeros.
\newblock {\em Stochastic Hybrid Systems}.
\newblock Taylor \& Francis, 2007.

\bibitem{TCP2}
D.~Chafaï, F.~Malrieu, and K.~Paroux.
\newblock On the long time behavior of the tcp window size process.
\newblock {\em Stochastic Process. Appl.}, 120(8):1518--1534, 2010.

\bibitem{Col}
V.~Colizza, A.~Barrat, M.~Barthélémy, and A.~Vespignani.
\newblock The role of the airline transportation network in the prediction and
  predictability of global epidemics.
\newblock {\em PNAS}, 103(7):2015--2020, 2006.

\bibitem{Dai}
J.~G. Dai.
\newblock On positive harris recurrence of multiclass networks: A unified
  approach via fluid limits models.
\newblock {\em The Annals of Applied Probability}, 5(1):49--77, 1995.

\bibitem{Dav93}
M.H.A. Davis.
\newblock {\em Markov Models and Optimization}.
\newblock Springer, 1993.

\bibitem{Del}
S.~Delattre, N.~Fournier, and M.~Hoffmann.
\newblock Hawkes processes on large networks.
\newblock {\em The Annals of Applied Probability}, 26(1):216--261, 2016.

\bibitem{DMT95}
D.~Down, Sean~P. Meyn, and R.~L. Tweedie.
\newblock Exponential and uniform ergodicity of markov processes.
\newblock {\em The Annals of Probability}, 23:1671--1691, 1995.

\bibitem{Dua}
A.~Duarte, A.~Galves, E.~Löcherbach, and G.~Ost.
\newblock Estimating the interaction graph of stochastic neural dynamics.
\newblock {\em arXiv:1604.00419}, 2016.

\bibitem{Ost}
A.~Duarte and G.~Ost.
\newblock A model for neural activity in the absence of external stimuli.
\newblock {\em Markov Processes and Related Fields}, 22(1):37--52, 2016.

\bibitem{Dum}
V.~Dumas, F.~Guillemin, and P.~Robert.
\newblock A markovian analysis of additive-increase multiplicative-decrease
  (aimd) algorithms.
\newblock {\em Advances in Applied Probability}, 2002.

\bibitem{Dut}
B.L. Dutta, P.~Ezanno, and E.~Vergu.
\newblock Characteristics of the spatio-temporal network of cattle movements in
  france over a 5-year period.
\newblock {\em Preventive Veterinary Medicine(1):79-94}, 2014.

\bibitem{Har}
J.~M. Harrison and S.~I. Resnick.
\newblock The stationary distribution and first exit probabilities of a storage
  process with general release rule.
\newblock {\em Mathematics of Operations Research}, 1:347--358, 1976.

\bibitem{Haw}
A.~Hawkes.
\newblock Spectra of some self-exciting and mutually exciting point processes.
\newblock {\em Biometrika}, 58:83--90, 1971.

\bibitem{Hon}
D.~Hong and D.~Lebedev.
\newblock Many tcp user asymptotic analysis of the aimd model.
\newblock {\em Rapport de recherche INRIA}, 4229, 2001.

\bibitem{Hos}
P.~Hoscheit, S.~Geeraert, H.~Monod, C.A. Gilligan, J.~Filipe, E.~Vergu, and
  M.~Moslonka-Lefebvre.
\newblock Dynamical network models for cattle trade: Towards economy-based
  epidemic risk assessment.
\newblock {\em Journal of Complex Networks}, 2016.

\bibitem{Jack}
J.R. Jackson.
\newblock Networks of waiting lines.
\newblock {\em Operations Research 5 (4): pp.518-521}, 1957.

\bibitem{Ker}
S.~Kernéis, R.~F. Grais, P.-Y. Boëlle, A.~Flahault, and E.~Vergu.
\newblock Does the effectiveness of control measures depend on the influenza
  pandemic profile?
\newblock {\em PLoS ONE}, 3(1), 2008.

\bibitem{Loc}
E.~Löcherbach.
\newblock Absolute continuity of the invariant measure in piecewise
  deterministic markov processes having degenerate jumps.
\newblock {\em arxiv.org/abs/1601.07123}, 2016.

\bibitem{Lev}
R.~Levins.
\newblock Some demographic and genetic consequences of environmental
  heterogeneity for biological control.
\newblock {\em Bulletin of the Entomological Society of America}, 15:237--240,
  1969.

\bibitem{Lib}
D.~Liberzon and A.~S. Morse.
\newblock Basic problems in stability and design of switched systems.
\newblock {\em IEEE Control Systems Magazine}, 19:59--70, 2001.

\bibitem{MAW}
E.~O. MacArthur, R. H. \&~Wilson.
\newblock {\em The theory of island biogeography}.
\newblock Princeton University Press, 2016.

\bibitem{Mal}
R.~Malhamé.
\newblock A jump-driven markovian electric load model.
\newblock {\em Advances in Applied Probability}, 22:564--586, 1990.

\bibitem{Mjack}
S.~P. Meyn and D.~Down.
\newblock Stability of generalized jackson networks.
\newblock {\em The Annals of Applied Probability}, 4(1):124--148, 1994.

\bibitem{MT1}
S.~P. Meyn and R.~L. Tweedie.
\newblock Stability of markovian processes i: Criteria for discrete-time
  chains.
\newblock {\em Advances in Applied Probability}, 24:542--574, 1992.

\bibitem{MT2}
S.~P. Meyn and R.~L. Tweedie.
\newblock Stability of markovian processes ii: Continuous-time processes and
  sampled chains.
\newblock {\em Advances in Applied Probability}, 25:487--517, 1993.

\bibitem{MT3}
S.~P. Meyn and R.~L. Tweedie.
\newblock Stability of markovian processes iii: Foster-lyapunov criteria for
  continuous-time processes.
\newblock {\em Advances in Applied Probability}, 25:518--548, 1993.

\bibitem{MT4}
S.~P. Meyn and R.~L. Tweedie.
\newblock State-dependent criteria for convergence of markov chains.
\newblock {\em The Annals of Applied Probability}, 4(1):149--168, 1994.

\bibitem{NGT}
Genki~I. Nagatani~T. and Tainaka K.
\newblock Epidemics of random walkers in metapopulation model for complete,
  cycle, and star graphs.
\newblock {\em Journal of Theoretical Biology}, 2018.

\bibitem{Pul}
H.~R. Pulliam.
\newblock Sources, sinks and population regulation.
\newblock {\em The American Naturalist}, 132(5):652--661, 1988.

\bibitem{Rit}
M.E. Ritchie.
\newblock {\em Wildlife and Landscape Ecology}, chapter Populations in a
  Landscape Context: Sources, Sinks, and Metapopulations, pages 160--184.
\newblock Springer, New York, 1997.

\bibitem{Verb}
J.~Verboom, K.~Lankester, and J.~A. Metz.
\newblock Linking local and regional dynamics in stochastic metapopulation
  models.
\newblock {\em Biological Journal of the Linnean Society}, 42(1-2):39--55,
  1991.

\bibitem{Wal}
J.~Walrand and P.~Varaiya.
\newblock Sojourn times and the overtaking condition in jacksonian networks.
\newblock {\em Advances in Applied Probability}, 12(4):1000--1018, 1980.

\end{thebibliography}
